\magnification 1000
\input amstex
\documentstyle{amsppt}
\vsize 8.25in
\voffset 1cm
\hoffset 1.3cm
\hoffset 4.5truemm
\topmatter

\rightheadtext{On partial algebraicity of holomorphic mappings}
\leftheadtext{Jo\"el Merker}
\title 
On the partial algebraicity of holomorphic\\ 
mappings between two real algebraic sets in \\
complex euclidean spaces of 
arbitrary dimension
\endtitle
\author Jo\"el Merker
\endauthor
\address Laboratoire d'Analyse, Topologie et Probabilit\'es,
Centre de Math\'ematiques et d'Infor\-ma\-tique, UMR 6632, 39 rue Joliot
Curie, F-13453 Marseille Cedex 13, France. Tel\,: 00 33 (0)4 91 11 35 50
\ Fax\,: 00 33 (0)4 91 11 35
52\endaddress
\email merker\@cmi.univ-mrs.fr
00 33 / (0)4 91 11 36 72 / (0)4 91 53 99 05
\endemail


\loadeufm

\def\M{{\Bbb M}}
\def\X{{\Bbb X}}
\def\L{{\Bbb L}} 

\def\N{{\Bbb N}}

\def\R{{\Bbb R}} 
\def\C{{\Bbb C}} 
\def\L{{\Bbb L}}

\def\v{\vert}

\def\dim{\text{\rm dim}}
\def\codim{\text{\rm codim}}
\def\1{{\text{\bf 1}}}

\abstract
The rigidity properties of the local invariants of real algebraic
Cauchy-Riemann structures imposes upon holomorphic mappings some
global rational properties (Poin\-car\'e 1907) or more generally
algebraic ones (Webster 1977). Our principal goal will be to unify the
classical or recent results in the subject, building on a study of the
transcendence degree, to discuss also the usual assumption of
minimality in the sense of Tumanov, in arbitrary dimension, without
rank assumption and for holomorphic mappings between two arbitrary
real algebraic sets.

\medskip
\noindent
{\smc R\'esum\'e.}  La rigidit\'e des invariants locaux des structures
de Cauchy-Riemann r\'eelles alg\'ebriques impose aux applications
holomorphes des propri\'et\'es globales de rationalit\'e (Poincar\'e
1907), ou plus g\'en\'eralement d'alg\'e\-bricit\'e (Webster 1977).
Notre objectif principal sera d'unifier les r\'esultats classiques ou
r\'ecents, gr\^ace \`a une \'etude du degr\'e de transcendance, de
discuter aussi l'hypoth\`ese habituelle de minimalit\'e au sens de Tumanov,
et ce en dimension quelconque, sans hypoth\`ese de rang et pour des
applications holomorphes quelconques entre deux ensembles
alg\'ebriques r\'eels arbitraires.
\endabstract 

\endtopmatter

\document

\head \S1. Introduction \endhead

The algebraicity or the rationality of local holomorphic mappings
between real algebraic CR manifolds can be considered to be one of the
most remarkable phenomena in CR geometry. Introducing the
consideration of Segre varieties in the historical article [18],
Webster generalized the classical rationality properties of
self-mappings between three-dimensional spheres discovered by
Poincar\'e and later extended by Tanaka to arbitrary dimension.
Webster's theorem states that biholomorphisms between Levi
non-degenerate real algebraic hypersurfaces in $\C^n$ are
algebraic. Around the eighties, some authors studied proper
holomorphic mappings between spheres of different dimensions or
between pieces of strongly pseudoconvex real algebraic hypersurfaces,
notably Pelles, Alex\-ander, Fefferman, Pinchuk, Chern-Moser,
Diederich-Fornaess, Faran, Cima-Suffridge, Forstneric, Sukhov, and
others (complete references are provided in [2,5,8,14,16,18,19]).  In
the past decade, removing the equidimensionality condition in the
classical theorem of Webster, Sukhov for mappings between Levi
non-degenerate quadrics [16], Huang [8] for mappings between strongly
pseudoconvex hypersurfaces, and Sharipov-Sukhov [14] for mappings
between general Levi non-degenerate real algebraic CR manifolds have
exhibited various sufficient conditions for the algebraicity of a
general local holomorphic map $f\: M\to M'$ between two real algebraic
CR manifolds $M\subset \C^n$ and $M'\subset \C^{n'}$. A necessary and
sufficient condition, but with a rank condition on $f$ is provided in
[2]. Recently, using purely algebraic methods, Coupet-Meylan-Sukhov
([5]\,; {\it see} also [6]) have estimated the transcendence degree of
$f$ directly. Building on their work, we aim essentially to study the
algebraicity question in full generality ({\it cf.}~Problem~2.4) and to
unify the various approaches of [2,5,8,12,14,18,19]. Notably, we shall
state necessary and sufficient conditions for the algebraicity of $f$
without rank condition and we shall study the geometry of the minimality 
assumption thoroughly.

\head \S2. Presentation of the main result \endhead

\subhead 2.1.~Algebraicity of holomorphic mappings and their transcendence
degree\endsubhead Let $U\subset \C^n$ be a small nonempty open
polydisc. A holomorphic mapping $ f\: U \to \C^{n'}$, 
$f\in \Cal H(U, \C^{n'})$, is called {\it
algebraic} if its graph is contained in an irreducible $n$-dimensional
complex algebraic subset of $\C^n\times \C^{n'}$. Using classical
elimination theory, one can show that, equivalently, each of its
components $g:=f_1, \ldots, f_{n'}$ satisfies a nontrivial polynomial
equation $g^ra_r+\cdots+ a_0= 0$, the $a_j\in \C[z]$ being
polynomials. We recall that a set $\Sigma\subset U$ is called real
algebraic if it is given as the zero set in $U$ of a finite number of
{\it real algebraic polynomials} in $(z_1,\ldots,z_n,\bar
z_1,\ldots,\bar z_n)$.

Let us denote by ${\Cal A}=\C[z]$ the ring of complex polynomial
functions over $\C^n$ and by ${\Cal R}=\C(z)$ its quotient field
$\hbox{Fr} \, ({\Cal A})$. By ${\Cal R}(f_1,\ldots,f_{n'})$ we
understand the field generated by $f_1,\ldots,f_{n'}$ over ${\Cal R}$,
which is a subfield of the field of meromorphic
functions over $U$ and which identifies with the collection of
rational functions
$$
R(f_1,\ldots,f_{n'})=P(f_1,\ldots,f_{n'})/ Q(f_1,\ldots,f_{n'}), \ \
\ P,Q\in {\Cal R}[x_1,\ldots,x_{n'}], \ \ Q\neq 0.
\tag 2.2
$$

Following [5], the {\it transcendence degree}
$\nabla^{tr}(f)$ of the field ${\Cal
R}(f_1,\ldots,f_{n'})$ with respect to the field
$\Cal R$ provides an integer-valued invariant measuring
the lack of algebraicity of $f$. In particular, {\it
$\nabla^{tr}(f)$ is zero if and only if $f$ is algebraic.} Indeed, by
definition $\nabla^{tr}(f)$ coincides with the maximal cardinal number
$\kappa'$ of a subset $\{f_{j_1},\ldots,f_{j_{\kappa'}}\}\subset
\{f_1,\ldots,f_{n'}\}$, $1\leq j_1 < \cdots < j_{\kappa'} \leq n'$
which is algebraically independent over ${\Cal R}$. In other words,
$\nabla^{tr}(f)=\kappa'$ means that there exists a subset
$\{f_{j_1},\ldots,f_{j_{\kappa'}}\}\subset
\{f_1,\ldots,f_{n'}\}$ such that there does not exist a nontrivial
relation $P(f_{j_1},\ldots,f_{j_{\kappa'}})\equiv 0$ in ${\Cal H}(U)$,
$P\in {\Cal R}[x_1,\ldots, x_{\kappa'}]\backslash \{0\}$, but that for
every $\lambda$, $\kappa'+1\leq \lambda \leq n'$, every $1\leq
j_1<\cdots<j_\lambda \leq n'$, there exists an algebraic relation
$Q(f_{j_1},\ldots,f_{j_{\lambda}})\equiv 0$, $Q\in {\Cal
R}[x_1,\ldots,x_{\lambda}]\backslash \{0\}$. Of course,
$\nabla^{tr}(f_1,\ldots,f_{n'})\leq n'$. An equivalent geometric
characterization of $\nabla^{tr}(f)$ states that
$\nabla^{tr}(f)=\kappa'$ if and only if the {\it dimension} of the
minimal for inclusion {\it complex algebraic} set $\Lambda_f\subset U
\times \C^{n'}$ containing the {\it graph} $\Gamma_f= \{(z,f(z))\in
U\times \C^{n'}\: z\in U\}$ of $f$ is equal to $n+\kappa'$ (this
complex algebraic set $\Lambda_f$ is of course necessarily
irreducible). In other words, $\nabla^{tr}(f)$ is an invariant
intrinsically attached to $f$ which is given with $f$ and which
possesses an algebraic {\it and} a geometric signification.  In a
metaphoric sense, we can think that $\nabla^{tr}(f)$ measures the {\it
lack of algebraicity} of $f$, or conversely, that it provides an
estimation of the {\it maximal partial algebraicity properties} of
$f$.

\subhead 2.3.~Presentation of the main result \endsubhead
Then, because the transcendence degree is an appropriate invariant,
more general than the dichotomy between algebraic and non-algebraic
objects, we shall as in [5] study directly the transcendence degree of
holomorphic mappings between two real algebraic sets.  Our main goal
is to provide a synthesis of the results in [2,5,8,14,18,19].  Thus,
quite generally, let $f\: U \to \C^{n'}$ be a local holomorphic
mapping sending an arbitrary irreducible real algebraic set $\Sigma
\cap U$ into another real algebraic set $\Sigma'\subset \C^{n'}$. As
the algebraicity of $f$ is a non-local property, we shall assume in
the sequel that $\Sigma\cap U$ is a smooth closed CR-submanifold of
$U$ and that there exists a second polydisc $U'\subset \C^{n'}$ with
$f(U)\subset U'$ such that $\Sigma'\cap U'$ is also a smoooth closed
CR-submanifold of $U'$. We shall denote by $M$ and $M'$ these
connected local CR pieces of $\Sigma$ in $U$ and of $\Sigma'$ in
$U'$. Of course, after shrinking again $U$ and $U'$, we can suppose
that $f$ is of constant rank over $U$. The topic of this article is to
study in full generality the following problem by seeking an optimal
bound\,:

\proclaim{Problem 2.4} 
Estimate $\nabla^{tr} (f)$ in terms of geometric invariants of $f$,
$\Sigma$, $\Sigma'$.
\endproclaim

\noindent
To begin with, we shall first assume that $M$ is somewhere
minimal in the sense of Tumanov, as in [2,5,8,\-14,\-18,19]. 
In the sequel, we shall say that a
property ${\Cal P}$ {\it holds at a Zariski-generic point $p\in \Sigma$} if
there exist a proper real algebraic subset $E$ of $\Sigma$, depending
on $\Cal P$, such that the property ${\Cal P}$ holds at each point of
$\Sigma\backslash E$. Let $\Delta$ be the unit disc in $\C$. Our main
result lies in the following statement from which we shall recover every
algebraicity result of the cited literature.

\proclaim{Theorem 2.5}
Assume that $M$ is CR-generic, connected and minimal in the sense of
Tumanov at a Zariski-generic point, and let $\Sigma''$ be the minimal
$($for inclusion, hence irreducible$)$ real algebraic set satisfying $f(M)
\subset \Sigma'' \subset \Sigma'$. Let $\kappa'$ denote the transcendence
degree $\nabla^{tr}(f)$ of $f$.  Then near a Zariski-generic point
$p''\in\Sigma''$, there exists a local algebraic coordinate system in
which $\Sigma''$ is of the form $\Delta^{\kappa'}\times
\underline{\Sigma}''$ for some real algebraic variety
$\underline{\Sigma}'' \subset \C^{n'-\kappa'}$.
\endproclaim

This theorem says that the degree of nonalgebraicity of $f$ imposes
some {\it degeneracy condition} on $\Sigma'$, namely to contain a
smaller real algebraic set $\Sigma''$ which is ``degenerate'' in the
sense that it can be locally straightened to be a product by a
complex $\nabla^{tr}(f)$-dimensional polydisc at almost every
point. The main interest of Theorem~2.5 lies in fact in its various
reciprocal forms which are listed in \S3 below. Of course, the
assumption that $\nabla^{tr}(f)$ equals an integer $\kappa'$ is no
assumption at all, since $\nabla^{tr}(f)$ is automatically given with
$f$, but in truth abstractly, namely in a non-constructive way, as is
$\Sigma''$. The only unjustified assumption with respect to Problem~2.4
is the minimality of $M$ in the sense of Tumanov and it remains also to
explain why the estimate given by Theorem~2.5 is sharp and
satisfactory.

Thus, let us firstly consider the sharpness.  If $\Sigma''$ is an
irreducible real algebraic set as above, it can be shown that there
exists the largest integer $\kappa_{\Sigma''}$ such that $\Sigma''$ is
a product of the form $\Delta^{\kappa_{\Sigma''}}\times
\underline{\Sigma}''$ near a Zariski-generic point in suitable
local algebraic coordinates ({\it see} Theorem~9.10). 
This integer will also be
abbreviated by $\kappa_{\Sigma''}$ and we shall say that $\Sigma''$ is
{\it $\kappa_{\Sigma''}$-algebraically degenerate}. We shall also
write $\Sigma''\cap V''\cong_{\Cal A}\Delta^{\kappa_{\Sigma''}}
\times \underline{\Sigma}''$ to mean that $\Sigma''$
intersected with the small open set $V''$ is equivalent to a product
in {\it complex algebraic} (abbreviation\,: $\Cal A$)
coordinates. Then Theorem~2.5 states in summary that
$\kappa_{\Sigma''} \geq \nabla^{tr}(f)$.  With this notion defined and
this rephrasing of Theorem~2.5 at hand, we now notice that it can of
course well happen that $\kappa_{\Sigma''} >
\nabla^{tr}(f)$. For instance, this happens when $n=n'$, when $f$ is
an algebraic biholomorphic map, so $\nabla^{tr}(f)=0$, and when
$\Sigma=\Sigma'=\Sigma''=\C^n$ simply. So what~? In case where
$\kappa_{\Sigma''} > \nabla^{tr}(f)$, by an elementary
observation we shall show that {\it a suitable perturbation of $f$ can
raise and maximize its possible transcendence degree}. The precise
statement, which establishes the desired sharpness, is as follows.

\proclaim{Theorem~2.6}
Assume that $f$ is nonconstant and that $\Sigma''$ is the {\rm minimal
for inclusion} real algebraic set satisfying $f(M)\subset
\Sigma''\subset \Sigma'$. Remember that $\Sigma''$ is locally equivalent
to $\Delta^{\kappa_{\Sigma''}}\times
\underline{\Sigma}''$ at a Zariski-generic point. Then there
exist a point $p\in M$ such that $f(p)=:p''\in\Sigma''$ is a
Zariski-generic point of $\Sigma''$ and a local holomorphic self-map
$\phi$ of $\Sigma''$ fixing $p''$ which is arbitrarily close to the
identity map, such that $\nabla^{tr}(\phi\circ f)=\kappa_{\Sigma''}$
(of course, because of Theorem~2.5, it is impossible to produce
$\nabla^{tr}(\phi\circ f)>\kappa_{\Sigma''}$).
\endproclaim

Secondly, let us discuss the (until now still unjustified) minimality
in the sense of Tumanov assumption. Remember that CR manifold without
infinitesimal CR automorphisms are quite exceptional and in any case
carry few self-maps. In \S13.2 
below, we shall observe the following.

\proclaim{Theorem~2.7}
Let $M$ be a nowhere minimal real algebraic CR-generic manifold
and assume that the space of infinitesimal CR-automorphisms of $M$ 
is nonzero. Then $M$ admits a local one-parameter family of nonalgebraic
biholomorphic self-maps.
\endproclaim

\head \S3. Five corollaries\endhead

The direct converse of the main Theorem~2.5 bounds $\nabla^{tr}(f)$ as
follows and gives an optimal sufficient condition for $f$ to be
algebraic.

\proclaim{Theorem~3.1} 
Let $f\in {\Cal H}(U,\C^{n'})$ with $f(M) \subset \Sigma'$ and assume
$M$ is CR-generic and minimal at a Zariski-generic point. Then
$\nabla^{tr}(f) \leq \kappa_{\Sigma''}=$ algebraic degeneracy degree
of the minimal for inclusion real algebraic set $\Sigma''$ such that
$f(M) \subset \Sigma'' \subset \Sigma'$.  In particular, $f$ is
necessarily algebraic if there does not exist an $1$-algebraically
degenerate real algebraic set $\Sigma'''$ with $f(M)\subset
\Sigma'''\subset \Sigma'$.
\endproclaim

As Theorem~2.6 shows that, otherwise, $f$ can be slightly perturbed to
be nonalgebraic, this theorem provides a {\it necessary and
sufficient} condition for $f$ to be algebraic. Thanks to this
synthetic general converse, we will recover results of the cited
litterature as corollaries.  We also obtain as a corollary the
celebrated algebraicity result in [18] from the following statement.

\proclaim{Corollary 3.2}
\text{\rm ([2])}
Any biholomorphism between pieces of smooth real algebraic CR-generic
holomorphically nondegenerate manifolds in $\C^n$ which are minimal at
a Zariski-generic point must be algebraic.
\endproclaim

\demo{Proof} 
Let $f\: M\to M'$ biholomorphic. Suppose by contradiction that
$\nabla^{tr}(f)=\kappa'\geq 1$. We prove in this paper (Corollary~9.14)
that $M'$ is holomorphically degenerate if and only if $M'$ is
algebraically degenerate. Since $f(M)\equiv M'$ is locally onto,
necessarily $M''=M'$ is $\kappa'$-holomorphically degenerate, by
Theorem~2.5, a contradiction.
\qed
\enddemo

For instance, here is another more general consequence
of Theorem~3.1. By $\Cal V_{\C^n}(p)$, we denote a small
open polydisc centered at $p$ whose size may shrink.
Let $\Sigma_{reg}'$ denote the regular part of 
$\Sigma'$ in the sense of real algebraic geometry.

\proclaim{Corollary~3.3} 
Let $f\in {\Cal H}(U,\C^{n'})$ with $f(M) \subset \Sigma'$ and assume
$M$ is CR-generic and minimal at a Zariski-generic point. If
$f(M\cap \Cal V_{\C^n}(p))) \supset {\Cal V}_{\C^{n'}}(f(p)\cap
\Sigma_{reg}')$ for some point $p\in M$, then 
$\Sigma''=\Sigma'$ and $\nabla^{tr}(f) \leq \kappa_{\Sigma'}$.
\endproclaim

Recall now that a CR-generic manifold $M$ is called {\it Segre-transversal}
at $p\in M$ if $\forall \ V= {\Cal V}_{\C^n}(p)$, $\exists \, k\in
\N_*$, $\exists \, q_1,\ldots,q_k\in V\cap S_{\bar{p}}$ such that
$T_pS_{\bar{q}_1}+\cdots+T_pS_{\bar{q}_k}=T_p\C^n$, the $S_{\bar{r}}$
denoting Segre varieties, hence $p\in S_{\bar{q}_1},\ldots,p\in
S_{\bar{q}_k}$. Segre-transversal $M$'s are always minimal. In fact,
Segre-transversality at $p\in M$ appears to be equivalent to
minimality at $p$ ({\it see} [5]) in codimension one, when $M$ is a
hypersurface, but in codimension $\geq 3$ there exist already some
minimal and not Segre-transversal $(M,p)$'s.  Since the question
whether minimal CR-generic manifolds $(M,p)$ of codimension two are
Segre transversal was left open in [5], we devote \S10 to answer it
affirmatively. Now, another important and direct consequence of our
main Theorem~2.5 is the following statement, implying in particular
[8].

\proclaim{Corollary~3.4}
\text{\rm ([5])}
Let $f\in {\Cal H}(U, \C^{n'})$ with $f(M)
\subset \Sigma'$ and assume that $M$ is
Segre-transversal at $p$. Then $\nabla^{tr}(f) \leq m':=$ the maximal
dimension of a complex algebraic variety $A_{f(q)}'$ with $f(q) \in
A_{f(q)}' \subset \Sigma'$, where $q$ runs in $M$.
\endproclaim

\demo{Proof} 
First, $M$ is minimal at $p$ ([5]; {\it see} also \S14 here). Clearly,
if $\kappa'=\nabla^{tr}(f)$, then $\Sigma'$ contains a piece of
algebraic set $A'\cong_{\Cal A} \Delta^{\kappa'}$ by Theorem~2.5, so
$\kappa' \leq m'$.
\qed \enddemo

A fourth consequence of Theorem~3.1 is\,:

\proclaim{Corollary~3.5} 
\text{\rm ([19])} 
Let $f\in {\Cal H}(U,\C^{n'})$ with $f(M \cap U) \subset \Sigma'$ and
assume that $M$ is minimal at a Zariski-generic
point. If $\Sigma'$ does not any contain complex
algebraic varieties, then $\nabla^{tr}(f)=0$, i.e. $f$ is algebraic.
\endproclaim

\demo{Proof} 
As above, if $\kappa'=\nabla^{tr}(f)\geq 1$, then $M'\supset
A'\cong_{\Cal A} \Delta^{\kappa'}$.
\qed
\enddemo

In conclusion, we therefore unify the
results in the papers [2,5,8,14,18,19].

\subhead 3.6.~Second reflection\endsubhead
Let $f\: M\to M'$, $M\: \rho(z,\bar{z})=0$, $0\in M$, $M'\:
\rho'(z', \bar{z}')= 0$, $0\in M'$, $\rho$, $\rho'$ real algebraic 
polynomials\,: $\rho\in \C[z,\bar{z}]^d$, $\rho'\in \C[z',
\bar{z}']^{d'}$, $d=\text{\rm codim}_{\R} M$, $d'= \text{\rm
codim}_{\R} M'$, and let $S_{\bar{w}}$, $S_{\bar{w}'}'$ denote Segre
varieties, $w\in {\Cal V}_{\C^n}(0):=U$, $w'\in {\Cal
V}_{\C^{n'}}(0):=U'$, $S_{\bar{w}}:=\{z\in U\: \rho(z,\bar{w})=0\}$,
$S_{\bar{w}'}':=\{z'\in U' \: \rho'(z', \bar{w}')=0\}$. Following [7,10,19],
for every subset $E' \subset U'$, we define the {\it first} and the
{\it second reflection} of $E'$ (which appears in an article
of Diederich-Fornaess, Ann. Math. {\bf 107} (1978), 371--384) by
$$
r_{M'}(E'):= \{w'\in U'\: S_{\bar{w}'}' \supset E'\}, \ \ \ \ \
r_{M'}^2(E') := r_{M'} (r_{M'}(E')).
\tag 3.7
$$ 
One establishes easily that given $f\: M\to M'$, then ({\it see}
[10])\,:
$$
f(z)\in \X_{z,\bar{w}}':= r_{M'} (f(S_{\bar{z}})) \cap
r_{M'}^2(f(S_{\bar{w}})).
\tag 3.8
$$
Although the paper [19] contains a slightly different statement using
the double reflection determination~\thetag{3.8}, we shall summarize
its main theorem by the following statement\,:
$$
\left< \,
\text{\rm dim}_\C \, \X_{z,\bar{w}}' = 0 \ \ \ \forall \ z, w \in {\Cal
V}_{\C^n}(0), \ \hbox{with} \ 
z\in S_{\bar{w}} \, \right> \ \Rightarrow \ \nabla^{tr}(f) = 0.
\tag 3.9
$$
In fact, to review the version [10], the condition $\text{\rm dim}_\C
\, \X_{z, \bar{w}}' = 0$ $\forall \ z, w\in {\Cal V}_{\C^n}(0), z\in
S_{\bar{w}}$, yields easily that $f$ is algebraic on every Segre
variety $S_{\bar{w}}$ and it suffices afterwards to apply
Theorem~5.1 below to get $\nabla^{tr}(f)= 0$. But the main point here
is again that, contrary to the one in Theorem~2.5, {\it this
condition~\thetag{3.9} is only sufficient}. Indeed, it appears that
there is no natural reason why there should exist a nonalgebraic
perturbation $\phi\circ f$ as in Theorem~2.6 in case where
$\dim_\C\X_{z,\bar w}'\geq 1$ for all $z,\, w, \, z\in S_{\bar w}$. 
For instance, there are cases where $\text{\rm dim}_\C
\X_{z,\bar{w}}' \geq 1$ $\forall \ z, w \in {\Cal V}_{\C^n}(0), z\in
S_{\bar{w}}$, but a suitable modification of the method in [7,10,19] yields
another determinacy set $\M_{z, \bar{w}}'$ such that $f(z)\in
\M_{z,\bar{w}}'$ again, such that $ \text{\rm dim}_\C
\M_{z,\bar{w}}'= 0$ $\forall \ z, w \in {\Cal V}_{\C^n}(0), z\in
S_{\bar{w}}$ and such that, moreover, $\M_{z,\bar{w}}'$ is also
appropriate for showing that $f$ is algebraic as in~\thetag{3.9}
above. By examples similar to the ones in [10], 
it can be shown that neither the
condition $\text{\rm dim}_\C \X_{z,\bar{w}}' \geq 1$ $\forall 
\ z, w \in {\Cal V}_{\C^n}(0), z\in S_{\bar{w}}$, nor $\text{\rm dim}_\C
\M_{z,\bar{w}}' \geq 1$ $\forall \ z, w \in {\Cal V}_{\C^n}(0), z\in
S_{\bar{w}}$ mean that we can find a nonalgebraic perturbation 
$\phi\circ f$\,: both conditions are {\it only sufficient}.
Further, the preprints [7,10] explains the unexpected
phenomena which are caused by the action of second reflection
$r_{M'}^2$.
By an examination of the examples in [10], the reader can
realize that {\it the trouble in them comes from the fact that $M'$ is not a
priori a piece of the minimal for inclusion real algebraic set
$\Sigma''$ containing $f(M)$}. And we can add that these puzzling phenomena
are essentially due to the fact that $r_{M'}$ {\it reverses inclusion of
sets}\,: $E'\subset F'$ implies $r_{M'} (E')\supset r_{M'} (F')$ and for
this reason, $E'\cap r_{M'} (E')$ and $F'\cap r_{M'}(F')$ cannot be
comparable {\it a priori}. For instance, if $f(M)\subset M''\subset M'$ are
as above, one has $r_{M'}(f(S_{\bar{z}}))\subset
r_{M''}(f(S_{\bar{z}}))$ but the two double reflection sets
$$
\X_{z,\bar{w}}':= r_{M'} (f(S_{\bar{z}})) \cap
r_{M'}^2(f(S_{\bar{w}})) \ \ \ \ \ \hbox{and} \ \ \ \ \
\X_{z,\bar{w}}'':= r_{M''} (f(S_{\bar{z}})) \cap
r_{M''}^2(f(S_{\bar{w}}))
\tag 3.10
$$
are {\it not comparable} in general. Explicit examples showing that one can
play very freely with these inclusions are given in [7,10]. This shows
that it is more natural to have a condition about $\X_{z,\bar{w}}''$
with $M''$ being a piece of the minimal for inclusion 
real algebraic set $\Sigma''$ which is {\it smooth}. We can of course 
assume smoothness after shrinking $M$ and $U$ a little
bit, since if otherwise $f(M)$ is contained in the (real algebraic)
singular part of $\Sigma''$ then $\Sigma''$ is not minimal 
for inclusion. Finally, after assuming that $M'$ is already
minimal for inclusion, we can derive quickly from our Theorem~2.5 the
contraposition of Theorem~1.1 in [19].

\proclaim{Corollary 3.11} 
Let $f\in {\Cal H}(U, \C^{n'})$, $f(M) \subset M'$, with $M$
CR-generic and minimal at a Zariski-generic point. Assume that $M'$ is
minimal for inclusion real algebraic CR-generic containing
$f(M)$ and smooth. If $\nabla^{tr}(f) \geq 1$, then there exists $D_M
\subset M\cap U$ a Zariski-open subset of $M$ such that $\forall \ p \in
D_M$, $\text{\rm dim}_\C \,\X_{z,\bar{w}}' \geq \nabla^{tr}(f) \geq
1$, $\forall \ z, w \in {\Cal V}_{\C^n}(p)$, $z\in S_{\bar{w}}$.
Equivalently, if  $\text{\rm dim}_\C \,\X_{z,\bar{w}}'=0$ on 
an open set, then $f$ is algebraic. In particular, if 
$\dim_\C \, r_{M'}(f(S_{\bar z}))=0$ on a open set, then 
$f$ is algebraic.
\endproclaim

\demo{Proof} Let $\kappa':= \nabla^{tr}(f)$. By Theorem~2.5, $M'$ is
at least $\kappa'$-algebraically degenerate, {\it i.e.} $\kappa_{M'}
\geq \kappa'$. Then $(M',p')\cong_{\Cal A}\Delta^{\kappa_{M'}}
\times \underline{M}'$ locally in a neighborhood of a Zariski-generic
point $p'\in M'$. Let $D_M$ be the open set of points $p\in M\cap U$
such that $M' \cap {\Cal V}_{\C^{n'}}(f(p))\cong_{{\Cal A}} 
\Delta^{\kappa_{M'}}\times
\underline{M}'$. Since $M'$ is minimal 
for inclusion, we claim that $D_M$ is Zariski-open in $M\cap U$. Indeed,
otherwise, $f(M\cap U)$ would be contained in the set of points $p'\in
M'$ where $M' \cap {\Cal V}_{\C^{n'}} (p')
\not\cong_{{\Cal A}} \Delta^{\kappa_{M'}}\times \underline{M}'$,
which is a {\it proper} real algebraic subvariety of $M'$ 
{\it see} Theorems~9.10 and~9.16 below, contradicting the choice of
$M'$. Assume therefore that $M'= \Delta^{\kappa_{M'}}\times
\underline{M}'$ in $U' ={\Cal V}_{\C^n}(0)$. 
Let $\pi'\: \C^{n'} \to \C^{n'-\kappa_{M'}} \times 0$ be the
projection. It is easy to show that for an arbitrary set $E' \subset
U'$, we have
$$
r_{M'}(E')=r_{\underline{M}'} (\Delta^{\kappa_{M'}}\times \pi'(E')), 
\ \ \ \ \ r_{M'}^2(E') = r_{\underline{M}'}
(r_{\underline{M}'} (\Delta^{\kappa_{M'}}\times \pi'(E'))).
\tag 3.12
$$
Consequently, $\X_{z,\bar{w}}' \supset\Delta^{\kappa_{M'}}\times \pi'(f(z))$ 
and $\text{\rm dim}_\C \, \X_{z,\bar{w}}' \geq
\kappa_{M'}\geq \kappa'= \nabla^{tr} (f) \geq 1$.
\qed
\enddemo

\subhead 3.12.~Summary of the proof of Theorem~2.5\endsubhead
Let $\nabla^{tr}(f)=:\kappa'$. Equivalently, the graph of
$f$ is contained in an algebraic $(n+\kappa')$-dimensional manifold in
$\C^n\times \C^{n'}$ and some $\kappa'$ components of $f$, say
$(f_1,\ldots,f_{\kappa'}):= f_{(\kappa')}$, make a transcendence basis
of the field extension $\Cal R\to \Cal R(f_1,\ldots, f_{n'})$. By the
assumption $f(M)\subset \Sigma'$, there exist algebraic relations
between $f$ and $\bar f$ and then between $f_{(\kappa')}$ and $\bar
f_{(\kappa')}$, after elimination. The main argument in
\S4 below shows that an algebraic dependence
$R'(z,\bar{z},f_{(\kappa')}(z),
\bar{f}_{(\kappa')}(\bar z)) \equiv 0$, for $z\in M$, 
can be transformed into an algebraic dependence
$S'(z,f_{(\kappa')})\equiv 0$ {\it which does not involve the
antiholomorphic components} in case $(M,p)$ is minimal. This fact
strongly relies on the Theorem 8.2 about propagation of partial
algebraicity along the Segre surfaces of $(M,p)$. Then $S'\equiv 0$,
and this will show easily that the set $\Sigma''$ is $\kappa'$-algebraically
degenerate.

\smallskip

\subhead 3.13.~Organization of the article\endsubhead 
\S4 provides the proof of Theorem~2.5, except the theorem on propagation of
algebraicity, to which \S5, \S6, \S7 and \S8 are devoted. \S7 contains
the proof of Theorem~2.7. \S9 presents the notion of algebraic
degeneracy. \S10 studies the notion of Segre-transversality in the
real algebraic context. Finally, \S11 shows how to produce a statement
equivalent to our main Theorem~2.5 using the so-called reflection
mapping.

\subhead 3.14.~Acknowledgement\endsubhead
Since I use transcendence degree in this paper, I would like to thank
Bernard Coupet, Francine Meylan and Alexander Sukhov, who kindly
provided me with their joint work in February 1998, and I wish to
address special thanks to Alexander for several interesting
discussions. Also, I have benefited of fruitful conversations about
the jet method with Bernard Coupet and Sylvain Damour. Finally, I
wish to thank the referee for his clever help.

\head \S4. Estimate of transcendence degree\endhead

\subhead 4.1.~Necessity in the main theorem\endsubhead
We proceed here first to the proof of the easy Theorem~2.6. Thanks to
the existence of algebraic stratifications and thanks to
delocalization ({\it i.e.}  choice of a smaller $U$ centered at
Zariski-generic point), we reduce the problem to a {\it nonconstant}
holomorphic map $f\: M\to \Delta^{\kappa'}\times
\underline{\Sigma}''$, $\kappa'\geq 1$.
To begin with, assume first that $f$ is algebraic, that $p=0$ and
$f(p)=0$. We write $f= (\underline{f}_1, \ldots,
\underline{f}_{n''}, f_1, \ldots, f_{\kappa'})$. 
Obviously, all perturbations $ (\underline{f}_1, \ldots,
\underline{f}_{n''}, g_1, \ldots, g_{\kappa'})$ of $f$, where the map
$g\: {\Cal V}_{\C^n}(0) \to
\Delta^{\kappa'}$, $g(0)= 0$, is an arbitrary holomorphic map, {\it
still send} $M$ {\it into} $\Delta^{\kappa'}\times
\underline{\Sigma}''$. We need to find such a perturbation of the form
$\phi\circ f$ with transcendence degree $\kappa'$.

Because $\Sigma''$ is assumed to be minimal for inclusion containing
$f(M)$, we have $f_1\not\equiv 0$. Next, we choose a transcendent
entire holomorphic function with high order of vanishing $\varpi\:
\Delta \to \Delta$, $\varpi(0)=0$, say for instance $\varpi(z_1)= 
(\sin z_1)^a$, $a\in \N_*$, $a >> 1$, and we define the map $\phi\: 
\Delta^{\kappa'}\times
\underline{\Sigma}''\to \Delta^{\kappa'}\times
\underline{\Sigma}''$ by
$\phi(\underline{z}_1,\ldots, \underline{z}_{n''}, z_1, \ldots,
z_{\kappa'}):=(\underline{z}_1,\ldots, \underline{z}_{n''},
z_1+\varpi(z_1),z_2+\varpi^{\circ 2}(z_1), \ldots, z_{\kappa'}+
\varpi^{\circ \kappa'} (z_1))$, where we denote $\varpi^{\circ k} := 
\varpi \circ \cdots \circ \varpi$ ($k$ times). Then we have 
$\nabla^{tr}(\varpi, \varpi^{\circ 2}, \ldots,
\varpi^{\circ k})= k$ for all $k\in\N_*$ and $\phi
\circ f= (\underline{f}_1, \ldots, \underline{f}_{n''}, f_1+
\varpi(f_1),\ldots, f_{\kappa'}+\varpi^{\circ \kappa'}(f_1))$. If $f$
was algebraic, then clearly $\nabla^{tr} (\phi\circ f)\geq \kappa'$
and the rank at $0\in \C^n$ (or the generic rank) of $f$ is preserved
after composition by $\phi$, because $\phi$ is arbitrarily close
to the identity map in a neighborhood of the origin.

Now, if $f$ was not algebraic, we choose instead functions
$\varpi_1,\ldots, \varpi_{\kappa'}\: \Delta \to \Delta$ with high order of
vanishing at $0$ such that $\nabla^{tr}(f_1+\varpi_1(f_1),
f_2+\varpi_2(f_1), \ldots, f_{\kappa'} + \varpi_{\kappa'}(f_1))= \kappa'$,
which is possible. Then again $\nabla^{tr}(\phi\circ f )
\geq \kappa'$. 
\qed 

\subhead 4.2.~Sufficiency in the main theorem\endsubhead 
We establish here Theorem~2.5. This paragraph will essentially follow
the lines of the main argument given in [5]. Let
$\kappa':=\nabla^{tr}(f)$ and choose a point $p\in M$ at which $M$ is
minimal such that $\Sigma''$ is smooth (in the sense of real algebraic
geometry) and CR at $f(p)$. This choice can be simply achieved by
avoiding some two proper real algebraic subvarieties in the source and
in the target. Of course, the real algebraic set $\Sigma''$ which is
minimal for inclusion satisfying $f(M) \subset
\Sigma'' \subset \Sigma'$ is irreducible and $f(M)$ is not
contained in $\Sigma_{sing}''$. Suppose we have proved that $\Sigma''$
is $\kappa'$-algebraically at $f(p)$. Then Theorem~9.10 below yields
that $\Sigma_{reg, CR}''$ is $\kappa'$-algebraically degenerate
everywhere. We are thus reduced to prove Theorem~2.6 at one point. To
summarize, thanks to the above simplifications, it is clear that it
suffices now to establish the following local statement.

\proclaim{Theorem~4.3} 
Assume that $f\: (M,p) \to (M',p')$ is a holomorphic map of constant
rank between two smooth CR-generic small manifold pieces of real
algebraic sets and let $\kappa':=\nabla^{tr}(f)$. If $(M,p)$ is
minimal and $(M',p')$ minimal for inclusion containing $f(M,p)$, then
$(M',p')$ is at least $\kappa'$-algebraically degenerate.
\endproclaim

\demo{Proof} 
For basic definitions about algebraic or transcendental extensions
that will be needed in this demonstration, we refer the reader to \S2
of [5] and the references therein. We only recall the following
important lemma. If $k$ is a subfield of a field $E$, then $E$ is said
to be an {\it extension field of} $k$. We write $k\to E$.

\proclaim{Lemma 4.4} 
Let $E=k(\alpha_1,\ldots,\alpha_n)$ be a \text{\rm finite}
extension of a field $k$, $n\in \N_*$, consisting of
rational functions of $(\alpha_1,\ldots,\alpha_n)$
with coefficients in $k$. Then $\nabla^{tr}(k\to
E)=\kappa'\geq 1$ if and only if, after renumbering\,:
\roster
\item"{\bf (1)}" 
$\alpha_1$ is transcendent over $k$, $\cdots$, $\alpha_{\kappa'}$ is
transcendent over $k(\alpha_1,\ldots,\alpha_{\kappa'-1})$\text{\rm ;}
\item"{\bf (2)}" 
$\alpha_{\kappa'+1},\ldots,\alpha_{n}$ are algebraic over
$k(\alpha_1,\ldots,\alpha_{\kappa'})$.
\endroster
\endproclaim

Let ${\Cal R}:=\text{\rm Fr}({\Cal A})=\C(z)$ denote the quotient field of
the ring of algebraic functions ${\Cal A}:=\C[z]$ over $\C^n$ and let
us consider ${\Cal R}(f_1,\ldots,f_{n'})$, the field of rational
functions generated by the components $f_1,\ldots,f_{n'}$ of $f$, {\it
i.e.} by fractions of the form
$$ 
\frac{P(f_1,\ldots,f_{n'})}{Q(f_1,\ldots,f_{n'})}, \
\ \ \ P(f_1,\ldots,f_{n'})= \sum_J a_J f^J, \ \ \ 
Q(f_1,\ldots,f_{n'})=\sum_J b_J f^J,
\tag 4.5
$$ 
$P$ and $Q$
being polynomials with coefficients $a_J\in {\Cal R}$, $b_J\in {\Cal
R}$. If $\kappa'$ denotes the transcendence degree of the field
extension ${\Cal R}
\to {\Cal R}(f)$, one can assume (after renumbering) that
$f_1,\ldots,f_{\kappa'}$ is the basis of transcendence, which means
({\it cf.}~4.4)

\roster
\item"{\bf 1.}" 
$f_1$ is transcendent over ${\Cal R}$, $\ldots$, $f_{\kappa'}$ is
transcendent over ${\Cal R}(f_1,\ldots,f_{\kappa'-1})$;
\item
"{\bf 2.}" 
$f_{\kappa'+1},\ldots,f_{n'}$ are algebraic over ${\Cal
R}(f_1,\ldots,f_{\kappa'})$.
\endroster

In particular {\bf 1} implies that there are no algebraic relations between
$(f_1,\ldots,f_{\kappa'})$. Let us denote 
$f_{(\kappa')}:=(f_1,\ldots,f_{\kappa'})$. Also,
{\bf 2} means that there are irreducible
monic polynomials
$$
S_{\kappa'+1}'(z,f_{(\kappa')}; X)\in {\Cal
R}(f_{(\kappa')})[X],\ldots,S_{n'}'(z,f_{(\kappa')}; X)\in {\Cal
R}(f_{(\kappa')})[X]
\tag 4.6
$$
such that writing
$$ 
S_j'=\sum_{0\leq k \leq m_j} S_{jk}'(z,f_{(\kappa')})X^{m_j-k}, \ \ \
\ S_{j0}'=1, \ \ \ \hbox{one has} \ \ \ S_j'(z,f_{(\kappa')};
f_j(z))\equiv 0
\tag 4.7
$$
identically for $z\in U$, $j= \kappa'+1,\ldots,n'$ with
$S_{jk}'(z,f_{(\kappa')})\in {\Cal R}(f_{(\kappa')})$. In other words,
the graph of $f$
$$
\Gamma(f)=\Gamma(f_1,\ldots,f_{n'})=\{(z,z')\in \C^n\times \C^{n'}; \
z_1'=f_1(z),\ldots,z_{n'}'=f_{n'}(z), z\in U\}
\tag 4.8
$$ 
is contained in the complex algebraic set
$$
\aligned
&
\Lambda'=\{(z,z')\in \C^n \times \C^{n'}; \ S_j'(z,z_{(\kappa')}';
z_j')=0,\\ & j=\kappa'+1,\ldots,n', z\in U, z'\in U'\} \subset \C^n
\times
\C^{n'},
\endaligned
\tag 4.9
$$ 
equipped with a natural projection $\tau'\: \C^n \times
\C^{n'} \to \C^n \times \C^{\kappa'}$, which is a local algebraic
biholomorphism $\tau'\: \Lambda' \backslash (\tau')^{-1}(\Upsilon') \to
\C^n \times \C^{\kappa'}$ outside the inverse image
$(\tau')^{-1}(\Upsilon')$ of the union $\Upsilon'$ of the discriminant
loci of the irreducible polynomials $S_j'$ (hence $\Upsilon'$ is a
complex algebraic subset of $\C^n \times \C^{\kappa'}$ of dimension
$\leq n+\kappa'-1$).

Consider now the graph of the transcendent basis
$$ 
\Gamma(f_1,\ldots,f_{\kappa'})= \{(z,z_{(\kappa')}')\in \C^n
\times \C^{\kappa'}; \ z_1'=f_1(z),\ldots,z_{\kappa'}'=f_{\kappa'}(z),
z\in U\}.
\tag 4.10
$$

\proclaim{Lemma 4.11} 
For any nonempty open set $V_M\subset M$, one has
$$
\Gamma(f_1,\ldots,f_{\kappa'})|_{V_M}
\not\subset \Upsilon'.
\tag 4.12
$$
\endproclaim

\demo{Proof} 
Assuming $V_M= V\cap M$ with both $V\subset \C^n$ and $V_M\subset M$ 
connected, it would imply $\Gamma(f_1,\ldots,f_{\kappa'})|_V
\subset \Upsilon'$ (identity principle), so there would be a
nontrivial algebraic relation between $f_1,\ldots,f_{\kappa'}$.
\qed
\enddemo

After delocalization, we then have $(p,f_{(\kappa')}(p))\not\in
\Upsilon'$ for all $p\in U$, so the complex algebraic
set $\Lambda'$ in $U\times U'$ can be locally
defined by equations of the form
$$
z_{\kappa'+1}'=h_{\kappa'+1}'(z,z_1',\ldots,z_{\kappa'}'),\ldots,
z_{n'}'=h_{n'}'(z,z_1',\ldots,z_{\kappa'}'),
\tag 4.13
$$ 
using the
algebraic implicit function theorem, where the 
$h_j'$ are {\it holomorphic algebraic} functions, {\it i.e.}
holomorphic functions whose graph is contained in a complex algebraic
set of dimension $n+\kappa'$.

By hypothesis,
$M'$ is given by real polynomial equations $P_j'(z',\bar{z}')=0, j =
1,\ldots,\sigma'$, near $p'$ and
we have $P_j'(f(z),\bar f(\bar z))\equiv 0$ for $z\in M$.
As in classical elimination theory, we
can replace the variables
$z_{\kappa'+1}',\bar{z}_{\kappa'+1}',\ldots,z_{n'}',\bar{z}_{n'}'$ by
the above values~\thetag{4.13} and their conjugate values (the $h_j'$
satisfy polynomial equations, so consider the other ``conjugate''
roots of these polynomials), take the product of these equations, use
Newton's identities and get that $\pi'(f(M))=f_{(\kappa')}(M)$ is
contained in a real algebraic set with polynomial equations 
$R_j'(z, \bar{z},
z_1',\ldots,z_{\kappa'}',\bar{z}_1',\ldots,\bar{z}_{\kappa'}')=0$,
$j=1,\ldots,\sigma'$.  We thus have
$$
R_j'(z,\bar z,f_{(\kappa')}(z),\bar f_{(\kappa')}(\bar
z))\equiv 0, \ \ \ j=1,\ldots, \sigma', \ z\in M.
\tag 4.14
$$

In summary, insofar we have eliminated the relatively algebraic
components $f_{\kappa'+1},\ldots,f_{n'}$ and we are now left with some
algebraic relations between the transcendence basis and its conjugate,
namely~\thetag{4.14} above. {\it The crucial Proposition~4.16 below
shows that all $R_j'\equiv 0$ necessarily}. We then claim that this
fact will readily imply that $M'$ contains (and is equal, by
minimality for inclusion) a $\kappa'$-algebraically degenerate set
like in Theorems~2.5 and 4.3.

Indeed, as $R_j'\equiv 0$, $1\leq j\leq \sigma'$, then for each $z\in
M$, the $\kappa'$-dimensional algebraic manifold
$A_z'=\{z_j'=h_j'(z,z_{(\kappa')}), \, j=\kappa'+1,\ldots,n'\}$ is
contained in $M'$. Let $\pi'\: (z,z')\mapsto z'$. 
Then the set of complex algebraic $A_z'$
parameterized by $z\in U$
$$ 
{\Cal C}=\{(z,z_{(\kappa')}',h_{(n'-\kappa')}'(z,z_{(\kappa')}'))\}
\subset U\times U'
\tag 4.15
$$ 
algebraically projects {\it via} \, $\pi'$ into $M'$, whenever $z\in
M$. The fibers of $\pi'|_{{\Cal C}}$ only depend on $z$, so there
exists an algebraic submanifold $N$ of $M$ where $\pi'$ has constant
rank and $\pi'\: {\Cal C} \cap (N \times U')\to M'$ is an algebraic
real diffeomorphism by minimality of $M'$ for inclusion.  Let
$Q_j(z,\bar z)=0$, $j=1,\ldots, \mu$ be some local real polynomial
equations for $N$. As $\Cal C\cap (N\times U')$ is given by the
equations $Q_j(z,\bar z)=0$,
$z_{(n'-\kappa')}'=h_{(n'-\kappa')}'(z,z_{(\kappa')}')$, the local
algebraic biholomorphism defined by $\tilde{z}_{(n'-\kappa')}:=
z_{(n'-\kappa')}'-h_{(n'-\kappa')}'(z,z_{(\kappa')}')$,
$\tilde{z}_{(\kappa')}':=z_{(\kappa')}'$, $\tilde{z}:=z$ clearly
straightens $\Cal C\cap (N\times U')$ to be the product of the real algebraic
manifold $\{(\tilde{z}, \tilde{z}_{(n'-\kappa')}')\:
Q_j(\tilde{z},\bar{\tilde{z}})=0, \, j=1,\ldots,\mu, \,
\tilde{z}_{(n'-\kappa')}'=0\}$ by the $\kappa'$-dimensional local 
polydisc $\{(0,\tilde{z}_{(\kappa')}',0)\}$. As $\pi'\: \Cal C\cap
(N\times U')\to M'$ is an algebraic CR-diffeomorphism, this shows that
$M'$ is at least $\kappa'$-algebraically degenerate. Granted
Proposition~4.16 below, then Theorem~4.3 will be proved.
\qed
\enddemo

It remains to show that no algebraic relation can be satisfied by the
transcendence basis $f_{(\kappa')}$ together with its conjugate $\bar
f_{(\kappa')}$. This is where the reflection principle and the
minimality of $M$ come on scene.

\proclaim{Proposition 4.16} 
Let $g:=(g_1,\ldots,g_{\kappa'})$, with $g_j(z)$ holomorphic in
$U$. Assume that $\nabla^{tr}(g)=\kappa'$. Then any polynomial
relation satisfied by $g$ and $\bar g$ over the field of rational
functions is trivial, namely $R(z,\bar{z}, g(z),\bar g(\bar z))\equiv
0$, $z\in M$, implies $R\equiv 0$.
\endproclaim

\demo{Proof} 
Let us proceed by contradiction. Then there exists an integer
$\mu\in\N_*$, such that we can write $R(z,\bar{z},g, \bar{g})=\sum_{1\leq i
\leq\mu} g^{\alpha_{i}} \ r_i(z,\bar{z},\bar g)$ where $\alpha_i\in
\N^{\kappa'}$ and where $r_i\neq 0$ for $i=1,\ldots,\mu$ and
such that $R(z,\bar z,g(z),\bar g(\bar z))\equiv 0$, for $z\in M$.
Here, $R$ and the $r_i$'s are polynomial. Of course, we can assume
that $\mu$ is the minimal integer for such a property. We shall use
the minimality of $\mu$ to derive the contradiction.  Here, since the
$r_i(z,\bar z,\bar g)$ are nonzero polynomials in $(z,\bar z,\bar g)$,
then for $z$ running in $U$, the terms $r_i(z,\bar z,\bar g(\bar z))$
(pull-back to $M$) can be considered, after obvious reordering, as
{\it polynomials in $z$ with coefficients being holomorphic functions
of $\bar z$} (we loose algebraicity in $\bar z$, because the terms
$\bar g(\bar z)$ are only holomorphic).  We can thus write them
as $r_i(\bar z;z):=\sum_J a_{i,J}(\bar z) \, z^J$, where such a sum is
understood to be finite. By minimality of $\mu$, necessarily no
$r_i(\bar z;z)$ vanishes identically on $M$.  For $i=2,\ldots,\mu$, we can
set $t_i(\bar z;z):=r_i(\bar z;z)/r_1(\bar z;z)$, which are terms of
the form $\sum_J a_{i,J}(\bar z) \, z^J/ (\sum_J a_{1,J}(\bar z) \,
z^J)$ and each sum is finite. 
Then we have the following relation for $z$ running over the
Zariski open subset $D_M:=\{r_1\neq 0\}$ of $M$\,:
$$ 
g^{\alpha_1}(z)+\sum_{2\leq i \leq \mu} t_i(\bar z; z) \
g^{\alpha_{i}}(z)\equiv 0, \ \ \ z\in D_M.
\tag 4.17
$$
Now, let $\bar L_1,\ldots,\bar L_m$ be a basis of $T^{0,1}M$ with
polynomial coefficients, where $m=\dim_{CR} M$. Of course, since
$g(z)$ is holomorphic, $\bar L_k(g_i)\equiv 0$ on $M$, for all $1\leq
k\leq m$ and $1\leq i\leq \kappa'$. Applying these CR derivations
to~\thetag{4.17}, we then see that the term $g^{\alpha_1}(z)$ is
automatically killed.  We therefore come to a dichotomy. Either $\bar
L_k(t_i(\bar{z};z))\equiv 0$, for all $1\leq k\leq m$ and all $2\leq
i\leq \mu$, $z\in D_M$, or there exists $1\leq k_*\leq m$ and $2\leq
i_*\leq \mu$ such that
$\bar L_{k_*}(t_{i_*}(\bar z ;z))\not\equiv 0$ on $M$. But this last
possibility would readily contradict the minimality of $\mu$.  Indeed,
as the coefficients of the $\bar L_k$'s are polynomial in $(z,\bar
z)$, all the terms $\bar L_{k_*}(t_{i}(\bar z ;z))$ are still of the
relatively rational form $\sum_J c_{i,J}(\bar z) \, z^J/ (\sum_J
d_{i,J}(\bar z) \, z^J)$ and after applying the CR derivation $\bar
L_{k_*}$ to~\thetag{4.17}, and after chasing the denominators, we would
obtain a similar polynomial relation $\sum_{2\leq i\leq
\mu} g^{\alpha_i}(z) \, r_i'(\bar z; z)\equiv 0$, $z\in M$,
$r_i'(\bar z,z)=\sum_J e_{i,J}(\bar z) \, z^J$, with a number
$\mu_1\leq \mu-1$ of terms strictly less than $\mu$, and this relation
is {\it nontrivial}, because $\bar L_{k_*}(t_{i_*}(\bar z
;z))\not\equiv 0$, which yields a contradiction in this case.  Thus,
we are left to discuss the first possibility, where $\bar
L_k(t_i(\bar{z};z))\equiv 0$, for all $1\leq k\leq m$ and all $2\leq
i\leq \mu$, $z\in D_M$. This case means that the real analytic
$t_i(\bar z;z)$'s are smooth and CR over $D_M$. Hence they admit a
holomorphic extension $s_i(z)$ to a neighborhood of $D_M$ in $\C^n$.
The important Proposition~4.18 below shows that this holomorphic
extension is in fact {\it holomorphic algebraic}.  But then
relation~\thetag{4.17} gives $g^{\alpha_1}(z)+\sum_{2\leq i \leq \mu}
s_i(z) \ g^{\alpha_{i}}(z)\equiv 0$, for $z\in D_M$, which is a
nontrivial algebraic relation between $(g_1,\ldots,g_{\kappa'})$, with
$\nabla^{tr}(g_1,\ldots,g_{\kappa'})=\kappa'$\,: this is again a
contradiction. Granted Proposition~4.18 below, then Proposition~4.16 
will be proved.
\qed 
\enddemo

\proclaim{Proposition 4.18} 
If $M$ is minimal at a Zariski-generic point and a relatively rational
with respect to $z$ function $t(\bar{z};z):=
\sum_J a_J(\bar z) \, z^J/(\sum_J b_J(\bar z) \, z^J)$, where the coefficients
$a_J(\bar z)$ and $b_J(\bar z)$ are holomorphic with respect to $\bar
z$, is CR on a Zariski open subset $D_M$ of $M$, then its holomorphic
extension $s(z)$ to a neighborhood of $D_M$ in $\C^n$ is algebraic.
\endproclaim

\demo{Proof} 
We localize first at a minimal point of $D_M$. Now assuming the
equation of $M$ is given by some local equations as in \S6.1 below, we
split the old coordinates $z$ in the new coordinates 
$(w,z)\in\C^m\times \C^d$ and we
complexify the equality $s(w,\bar z+i\bar{\Theta}(w,\bar w,\bar
z))\equiv \sum_{I,J} a_{I,J}(\bar w,\bar z) \, w^I z^J/(\sum_{I,J}
b_{I,J}(\bar w, \bar z) \, w^Iz^J)$, where each sum is finite, which
is valuable for $(w,z)\in D_M$, to obtain
$$
s(w,\xi+i\bar{\Theta}(w,\zeta,\xi))\equiv
\left[\sum_{I,J} a_{I,J}(\zeta,\xi) \, 
w^I z^J \left/\left(\sum_{I,J} b_{I,J} (\zeta, \xi) \,
w^Iz^J\right)\right.\right]_{z:=\xi+i\bar \Theta(w,\zeta,\xi)}
\tag 4.19
$$
identically on the extrinsic complexification ${\Cal M}$, {\it i.e.}
identically as power series in $(w,\zeta,\xi)$. Of course,
to complexify, we use the fact that $M={\Cal M} \cap \{\zeta =
\bar{w}, \xi= \bar{z}\}$ is maximally real in ${\Cal M}$, hence a
uniqueness set. Then~\thetag{4.19} shows that $s$ is algebraic on
each Segre surface ${\Cal S}_{\zeta_p,\xi_p}$ and $\underline{\Cal
S}_{w_p,z_p}$ of ${\Cal M}$, {\it see}~\thetag{6.4} for their
definition. We shall prove in Theorem~5.1 below that a holomorphic
function defined in a neighborhood of a minimal CR-generic manifold is
algebraic if and only if all its restrictions to Segre varieties are
algebraic. Therefore $s$ is algebraic, as desired. Granted Theorem~5.1
below, this will complete the proof of Theorem~4.3.
\qed 
\enddemo

\head \S5. Vector fields with complex algebraic flow \\
and partial algebraicity of holomorphic mappings \endhead

To conclude the proof of Theorem~4.3 above, it remains thus to
establish the following statement about separate algebraicity of
holomorphic mappings, which has been established by Sharipov-Sukhov
([14], {\it see} also [5] and \S10 below) in the case where $M$ is 
Segre-transversal, instead of being minimal.

\proclaim{Theorem 5.1} 
Let $g\in {\Cal H}({\Cal V}_{\C^n}(M), \C)$ be a holomorphic function,
$M$ being a real algebraic CR-generic manifold which is minimal at a
Zariski-generic point. Then $g$ is algebraic if and only if its
restriction to each Segre variety of $M$ is algebraic.
\endproclaim

After complexifying $M$ in a neighborhood of a minimal point, we shall
be naturally led to deduce Theorem~5.1 from a more general
statement. By slight abuse of terminology, we will call a
Frobenius-integrable $k$-dimensional distribution $L$ over a complex
manifold a ``{\it $k$-vector field}''. Such a distribution is said to
have ``{\it complex algebraic $k$-flow}'' if the local foliation that
it induces (by Frobenius' theorem) coincides with a trivial product
foliation by the $k$-dimensional polydiscs $\bigcup_{\theta}
\Delta^k\times \{\theta\}$ {\it in some local complex algebraic
coordinate system}. Its leaves will be called ``{\it k-curves}''.  We
refer the reader to \S7.1 for further material about $\L$-orbits.

\proclaim{Theorem~5.2} 
Let $\L=\{L^{\alpha}\}_{\alpha\in A}$ be a system of $k_\alpha$-vector
fields with complex algebraic flow on a small open connected set $U\subset
\C^n$. Then
\roster
\item"{\bf (1)}" 
For $p\in U$,
the $\L$-orbits $\Cal O_\L(U,p)$ are complex algebraic manifolds.
\item"{\bf (2)}"
A holomorphic function $g\in {\Cal H}(U, \C)$ is algebraic on each
$\L$-orbit if and only if it is algebraic on each $($complex
algebraic$)$ integral $k_\alpha$-curve of every element of $\L$.
\endroster
\endproclaim

In particular, if ${\Cal O}_{\L}(U, p)$ contains an open subset of
$\C^n$ for some $p\in U$, then a function $g\in {\Cal H}(U,\C)$ is
algebraic on $U$ under the sole assumption that it is algebraic on
$\L$-integral curves. This theorem, as well as its preliminary version
given by Sharipov-Sukhov, generalizes the well known separate
algebraicity principle in $\C^n$ proved in the book of Bochner-Martin
[4]\,: {\it A holomorphic function $g\in {\Cal H}(\Delta^n, \C)$ is
algebraic if and only if its restriction to every coordinate discs is
algebraic}.  As a corollary to this Theorem~5.2, we shall also provide
a new proof of the following theorem ({\it cf.} [2]).

\proclaim{Theorem 5.3}
The CR orbits of a real algebraic CR manifold are algebraic.
\endproclaim

We begin by explaining how Theorem~5.2 applies to provide a proof
of Theorem~5.1. For this, we follow and summarize the constructions of [9].

\head \S6. Foliations by complexified Segre varieties 
and algebraicity of CR-orbits of an algebraic CR manifold
\endhead

\subhead 6.1.~The extrinsic complexification of $M$\endsubhead
Let $M\subset \C^n$ be a real algebraic CR-generic submanifold and set
$m:=\dim_{CR} \, M$, $d:= \codim_\R \, M$, $m+d=n$.  Using the theory
of algebraic functions, one can see that there exist holomorphic
coordinates $t=:(w,z)\in\C^m\times \C^d$ and a {\it holomorphic
algebraic} $d$-vectorial function $Q(\bar w,t)=(Q_l(\bar w,t))_{1\leq
l\leq d}$ such that $M$ is given by the two equivalent systems of
cartesian defining functions $z=\bar Q(w,\bar t)$ or $\bar z=Q(\bar
w,t)$. As $M$ is real, these two systems of equations are equivalent
and there exists an invertible $d\times d$ matrix power series
$a(t,\bar t)$ such that $z-\bar Q(w,\bar t)\equiv a(t,\bar t) \, (\bar
z-Q(\bar w,t)$. We can furthermore assume that $T_0M=\C_w^m\times
\R_x^d$, in which case we shall write the equations of $M$ as
follows\,: $z=\bar z+i\bar\Theta(w,\bar w,\bar z)$ or $\bar
z=z-i\Theta(\bar w,w,z)$, where $\Theta$ vanishes to second order at
the origin. In substance, such functions are locally
holomorphic functions whose graph is contained in a complex algebraic
set of the same dimension as the basis. Let now $\tau:=(\bar
t)^c=:(\zeta,\xi)\in\C^m\times
\C^d$ denote the complexification variable of $\bar t$. Then the extrinsic
complexification of $M$ is given by the two, again equivalent, systems
of cartesian defining equations $z=\bar Q(w,\tau)$ or $\xi=Q(\zeta,t)$.
Following [9], we recall that there exist two systems of $m$-vector fields
$\Cal L=(\Cal L^1,\ldots,\Cal L^m)$ and $\underline{\Cal L}=
(\underline{\Cal L}^1,\ldots, \underline{\Cal L}^m)$, which are by 
definition the {\it complexifications} of two conjugate systems of 
$m$-vector fields spanning $T^{1,0}M$ and $T^{0,1}M$. These two systems
$\Cal L$ and $\underline{\Cal L}$ span two integrable subbundles
of $T\Cal M$ and they induce therefore two flow (Frobenius) foliations 
$\Cal F_{\Cal L}$ and $\Cal F_{\underline{\Cal L}}$ of $\Cal M$, thanks
to the integrability condition $[\Cal L^i,\Cal L^j]\subset \Cal L$
and $[\underline{\Cal L}^i,\underline{\Cal L}^j]\subset \underline{\Cal L}$.
As in [9], this fact can be rendered visible and straightforward
just by writing in coordinates these two $m$-vector fields, using
a vectorial and symbolic notation\,:
$$
\Cal M\: \ \ \
\Cal L =\frac{\partial }{\partial
w}+\bar Q_w (w,\zeta,\xi)
\ \frac{\partial }{\partial z} \ \ \ \
\ \hbox{and} \ \ \ \ \ \underline{\Cal L}=\frac{\partial }{\partial
\zeta}+Q_\zeta(\zeta,w,z) \frac{\partial }{\partial \xi}.
\tag 6.2
$$ 
Secondly, it is easy to observe that the two {\it different} (exercise) 
Segre varieties and conjugate Segre varieties as defined in [9], 
which can be rewritten as
$$
\left\{
\aligned
&
S_{\bar{t}_p}:=\{(w,z)\:
z=\bar{z_p}+i\bar{\Theta}(w,\bar{w}_p,\bar{z}_p)\} \ \
\ \ \ \hbox{and}\\
&
\overline{S}_{t_p}:=\{(\bar w,\bar z)\: 
\bar{z}=z_p-i\Theta(\bar w,w_p,z_p)\}
\endaligned\right.
\tag 6.3
$$ 
admit two different complexifications in $\Cal M$, which
we will denote by $\Cal S_{\tau_p}$ and
$\Cal S_{t_p}$, and which can be written as follows\,:
$$
\left\{
\aligned
&
\Cal S_{\tau_p}=
\Cal S_{\zeta_p,\xi_p}\: \ \ \ \zeta=\zeta_p, \ \xi
=\xi_p, \ z=\xi_p+i\bar{\Theta}(w,\zeta_p,\xi_p) \ \ \ \ \ \ \ \ \ \ 
\hbox{and} \\
&
\underline{\Cal S}_{t_p}=
\underline{\Cal S}_{w_p,z_p}\: \ \ \ w=w_p, \ z=z_p, \
\xi=z_p-i\Theta(\zeta,w_p,z_p).
\endaligned\right.
\tag 6.4
$$
Here, the coordinates $t_p$, $\bar t_p$ or $\tau_p$ are thought to be
{\it fixed} and the equations~\thetag{6.3} define two $m$-dimensional
complex submanifolds of $\Cal M$ which coincide in fact with the
coordinate intersections $\Cal M\cap \{\tau=\tau_p\}$ and $\Cal M\cap
\{t=t_p\}$ respectively. Applying $\Cal L$ to the equation of $\Cal
S_{\tau_p}$ and $\underline{\Cal L}$ to the equations of
$\underline{\Cal S}_{t_p}$, we see that $\Cal L$ is tangent to $\Cal
S_{\tau_p}$ and that $\underline{\Cal L}$ is tangent to
$\underline{\Cal S}_{t_p}$. We can thus summarize these
observations\,:
\roster
\item"{\bf 1.}"
The $\Cal S_{\tau_p}$ and the $\underline{\Cal S}_{t_p}$ form {\it
families of integral complex algebraic manifolds for $\Cal L$ and
$\underline{\Cal L}$} respectively.
\item"{\bf 2.}"
The $\Cal S_{\tau_p}$ are leaves of the flow foliation $\Cal F_{\Cal
L}$ of $\Cal M$ by $\Cal L$ and the $\underline{\Cal S}_{t_p}$ are
leaves of the flow foliation $\Cal F_{\underline{\Cal L}}$ of $\Cal M$
by $\underline{\Cal L}$.
\item"{\bf 3.}"
Therefore, the two $m$-flows of $\Cal L$ and of $\underline{\Cal L}$
are both complex algebraic, because the families $
\Cal F_{\Cal L}=\bigcup_{\tau_p\in\C^n, \v \tau_p\v <\delta}
{\Cal S}_{\tau_p}$ and $\Cal F_{\underline{\Cal
L}}=\bigcup_{t_p\in\C^n, \v t_p\v <\delta}
\underline{\Cal S}_{t_p}$ are clearly tangentially and
transversally algebraic.
\endroster
Finally, we will need the following straightforward lemma about
complexifications of Lie algebras including a well known
characterization of finite type (orbit-minimality) at a point,
valuable in the real analytic and real 
algebraic categories ({\it see} also [9]). Let 
$p^c:=(p,\bar p)\in\Cal M$ be the complexification of the
point $p$. Then\,:

\proclaim{Lemma~6.5}
The following five properties are equivalent\text{\rm \,:}
\roster
\item"{\bf (a)}" \
$\text{\rm Lie}_p \ (T^{1,0}M, T^{0,1} M)=\C\otimes T_pM$.
\item"{\bf (b)}" \
$\text{\rm Lie}_p \ (T^cM)=T_pM$.
\item"{\bf (c)}" \
$(M,p)$ is $T^cM$-orbit-minimal.
\item"{\bf (d)}" \
$\text{\rm Lie}_{p^c} \ (\Cal L, \underline{\Cal L})=T_{p^c}\Cal M$.
\item"{\bf (e)}" \
$(\Cal M,p^c)$ is $\{\Cal L, \underline{\Cal L}\}$-orbit-minimal. \qed
\endroster
\endproclaim

\subhead 6.6.~Deduction of Theorem~5.1 \endsubhead
As announced in \S5, Theorem~5.1 can be therefore deduced from the
more general statement Theorem~5.2, thanks to the following two
facts. Let there be given a holomorphic function $g\in\Cal H(\Cal
V_{\C^n}(p_0),\C)$ defined in a neighborhood of a minimal point $p_0$
of $M$, which we think to be the origin in the above coordinates
$(w,z)$ for $M$. Then\,:
\roster
\item"{\bf 1.}"
This function $g$ induces a function $g^c\: \Cal M\to \C$ defined by
$g^c(w,z,\zeta,\xi):=g(w,z)$. Since by assumption $g^c$ is algebraic
on every Segre variety $S_{\bar t_p}$, then $g^c$ is algebraic on
every complexified Segre variety $\Cal S_{\tau_p}$. On the other hand,
$g^c$ is also clearly algebraic on every conjugate complexified Segre
variety $\underline{\Cal S}_{t_p}$, because it is then {\it constant}
equal to $g(w_p,z_p)$.
\item"{\bf 2.}"
According to Lemma~6.5, $\Cal M$ is orbit-minimal at $(p_0,\bar p_0)$
for the system
$\L:=\{\Cal L,\underline{\Cal L}\}$ composed of (only) two $m$-vector
fields. Because by step {\bf 1} above, $g^c$ is algebraic on each
$m$-integral manifold of this system, then Theorem~5.2 clearly yields
the algebraicity of $g^c$ over $\Cal M$, whence $g$ is algebraic. \qed
\endroster

\subhead 6.7.~Deduction of Theorem~5.3 \endsubhead
Theorem~5.3 can also be deduced from Theorem~5.2, because of the
following relation between the CR-orbits $\Cal O_{CR}(M,p)=
\Cal O_{L,\bar L}(M,p)$ of points $p$ in $M$ and the $\{\Cal L,
\underline{\Cal L}\}$-orbits $\Cal O_{\Cal L,\underline{\Cal L}}
(\Cal M,p^c)$ of points $p^c\in\Cal M$, which is established in [9].
Let $\pi_t\: (t,\tau)\mapsto t$ denote the projection on the first
coordinate space, which is defined over the complexification space
$\C_t^n\times \C_\tau^n$. Let $\underline{\Lambda}:=\{(t,\tau)\:
\tau=\bar t\}$ denote the antidiagonal.

\proclaim{Proposition~6.8}
\text{\rm ([9])}
We have the following reciprocal complexication relations\,:
\roster
\item"{\bf (1)}" 
$\Cal O_{CR}(M,p)^c=\Cal O_{\Cal L, \underline{\Cal
L}}(\Cal M, p^c)$.
\item"{\bf (2)}" 
$\Cal O_{CR}(M,p)=\pi_t(\underline{\Lambda}\cap
\Cal O_{\Cal L, \underline{\Cal L}} (\Cal M, p^c))$.
\endroster
\endproclaim

Granted this proposition and thanks to the observation that the the
$m$-flows of $\Cal L$ and of $\underline{\Cal L}$ are naturally
algebraic, it now follows from Theorem~5.2 {\bf (1)} that the orbits
$\Cal O_{\Cal L, \underline{\Cal L}}(\Cal M, q)$ of various points 
$q\in \Cal M$ are all algebraic, whence with $q:=p^c$, we deduce from
{\bf (2)} above that the CR-orbits $\Cal O_{CR}(M,p)$ are all real
algebraic.  To conclude the proof of Theorem~2.5, it now only remains
to establish Theorem~5.2.

\head \S7. Algebraicity of orbits of vector fields having algebraic
flow and local holomorphic mappings of nowhere minimal CR-generic manifolds
\endhead

\subhead 7.1.~Definition of the concatenated flow maps \endsubhead
As in Theorem~5.2, let $\L=\{L^{\alpha}\}_{\alpha\in A}$ be a finite
set of nonzero holomorphic $k_\alpha$-vector fields defined on an open
connected set $U\subset \C^n$. For the moment's discussion, we do not
assume algebraicity of their flows. As all our reasonings will be
local, we shall assume that $U=\Delta^n$. Of course, after
demultiplying any such $k_\alpha$-vector in some $k_\alpha$ linearly
independent $1$-vector fields, we can also assume that all the
$L^\alpha$'s are (usual) $1$-vector fields. Thus, denote these vector
fields by $L^{\alpha}=\sum_{j=1}^n a_{j\alpha}\frac{\partial}{\partial
z_j}$, where $a_{j\alpha}(z)\in {\Cal H}(\Delta^n)$. By a well known
theorem, the global flow of each such vector field $L^\alpha\in \L$,
say $\phi_\alpha\: (t,p)\mapsto \exp \, (t \, L^\alpha)(p)\in
\Delta^n$, is defined over a certain maximal subdomain $\Omega_\alpha$ of 
$\C\times\Delta^n$ which contain $\Delta^n\times \{0\}$ 
and this global flow is a holomorphic map over that subdomain. 
However, following our 
version of Sussmann's constructions elaborated in [9] in
the (paradigmatic) analytic context, we shall proceed as follows
in order not to deal with these domains of the global flows
$\phi_\alpha$'s. Let us now denote by $(t,p)\mapsto L_t^\alpha(p)$
these flow maps in the sequel. To define $\L$-orbits, we shall need to
consider concatenated flow maps of the form $L_{t_k}^{\alpha_k}\circ
\cdots \circ L_{t_1}^{\alpha_1}(p)$. By definition, the $\L$-orbit
$\Cal O_\L(\Delta^n,p)$ of a point $p\in
\Delta^n$ is just the set of all such elements $L_{t_k}^{\alpha_k}\circ
\cdots \circ L_{t_1}^{\alpha_1}(p)\in\Delta^n$ where
$\alpha_1,\ldots,\alpha_k\in A$ and where $k\in\N_*$ is {\it a priori}
unbounded. Nevertheless, after comparison with the (semi-global)
constructions in [17], it appears that, due to the fact
that the $L^\alpha$'s are defined over $\Delta^n$, this set of points
gives nothing more interesting for $k>3n$. According to Nagano's
theorem (revisited by Sussmann with the above definitions), the main
property of $\L$-orbits lies in the fact that they are closed complex
analytic submanifolds of $\Delta^n$ passing through $p$ and for this,
$k\leq 3n$ suffices. Furthermore, there is another phenomenon which
is due to the principle of analytic continuation and which is quite
well known ({\it cf.} [17])\,: the local and the global CR-orbits
coincide locally in the real analytic category, which is false in the
smooth $\Cal C^\infty$ category. Therefore, to study orbits and
to fix the domains of the concatenated flow maps, it suffices to
choose some positive number $\delta>0$ such that for all $p\in {1\over
2^n}\Delta^n$, all $k\in\N_*$ with $k\leq 3n$ and all $t_1,\ldots,t_k$
with $\v t_j\v \leq \delta$, one has $L_{t_k}^{\alpha_k}\circ \cdots
\circ L_{t_1}^{\alpha_1}(p)\in\Delta^n$, which is clearly possible by
continuity of the flow maps. In summary, there is no restriction to
define from the beginning the $\L$-orbits of various points $p\in
{1\over 2^n}\Delta^n$ as the set of all elements $L_{t_k}^{\alpha_k}\circ
\cdots \circ L_{t_1}^{\alpha_1}(p)\in\Delta^n$ where
$\alpha_1,\ldots,\alpha_k\in A$, where $k\leq 3n$ and where $\v t_j\v
\leq \delta$. We shall say that the point $p$ is $\L$-orbit-minimal if
its $\L$-orbit (with this definition) contains an open neighborhood of
$p$ in $\Delta^n$.

\subhead 7.2.~Algebraicity of the $\L$-orbits\endsubhead
Assuming now that all the flows of the elements of $\L$ are complex
algebraic, we can easily check part {\bf (1)} of Theorem~5.2.
Indeed, according to the construction given in [9], for
an arbitrary fixed point $p\in {1\over 2^n}\Delta^n$, there exists
a integer $e_p\leq n$ such that
\roster
\item"{\bf 1.}"
The holomorphic maps $(t_1,\ldots,t_k)\mapsto 
L_{t_k}^{\alpha_k}\circ \cdots
\circ L_{t_1}^{\alpha_1}(p)\in\Delta^n$ are of generic rank
equal to $e_p$ for all $e_p\leq k\leq 3n$.
\item"{\bf 2.}"
There exist $t_1^*,\ldots,t_{e_p}^*\in\C$ arbitrarily close to $0$ such that
the map 
$$
\Gamma_{e_p}\: (t_1,\ldots,t_{e_p}) \mapsto
L_{-t_1^*}^{\alpha_1}\circ \cdots\circ L_{-t_{e_p}^*}^{\alpha_{e_p}}
\circ L_{t_{e_p}}^{\alpha_{e_p}}\circ \cdots \circ 
L_{t_1}^{\alpha_1}(p)\in\Delta^n
\tag 7.3
$$ 
is of {\it constant} rank $e_p$ over a neighborhood $\Cal T^*$ of 
$(t_1^*,\ldots,t_{e_p}^*)$ in $(\delta\Delta)^{e_p}$.
\item"{\bf 3.}"
The map $\Gamma_{e_p}$ (clearly)
satisfies $\Gamma_{e_p}(t_1^*,\ldots,t_{e_p}^*)=p$.
\item"{\bf 4.}"
The image $\Gamma_{e_p}(\Cal T^*)$ is an $e_p$-dimensional submanifold
of $\Delta^n$ through $p$ which coincides with the $\L$-orbit 
$\Cal O_\L(\Delta^n,p)$ of $p$ in a neighborhood of $p$ in $\Delta^n$.
\endroster

\noindent
Notice that this statement clearly shows that the $\L$-orbit of $p$ is 
a {\it submanifold} and that it provides this orbit with the regularity 
of the concatenated flow map $\Gamma_{e_p}$. 

\demo{End of proof of Theorem~5.2}
As we supposed that the flows are complex algebraic, the above maps
$\Gamma_{e_p}$ are algebraic and then clearly the orbits are complex
algebraic by {\bf 4}, {\it q.e.d.} Of course, similar other regularity
properties of orbits in other differentiable categories rely upon the
regularity $\Cal C^\infty$, $\Cal C^k$, $\Cal C^{l,\alpha}$ of the
flow maps. This proves part {\bf (1)} of Theorem~5.2 and we shall
establish part {\bf (2)} in \S8 below.
\qed
\enddemo

\subhead 7.4.~Flow-bow theorem \endsubhead
Now, as in the differentiable theory, it is easy to deduce from the
above four properties {\bf 1,2,3,4} a {\it complex algebraic} flow-box
theorem, that will be useful.

\proclaim{Theorem~7.5} 
Let $p\in U$ and set $e_p=\text{\rm dim} {\Cal O}_{\L}(\Delta^n,p)$.
Then there exist a neighborhood $V$ of $p$ in $\Delta^n$ and a
biholomorphism $\Phi\: V\to \Delta^{e_p} \times \Delta^{n-e_p}$ with complex
algebraic components such that
\roster
\item"{\bf (1)}"
$\Phi^{-1}(\Delta^{e_p} \times \{0\})= {\Cal O}_{\L} (\Delta^n,p) \cap V$.
\item"{\bf (2)}" $\forall \ \zeta_{(n-e_p)} \in \Delta^{n-e_p}$,
$\Phi^{-1}(\Delta^{e_p}\times \{\zeta_{(n-e_p)}\})$ is contained in a single
$\L$-orbit.
\item"{\bf (3)}" The function $\Delta^n\ni p\mapsto
\text{\rm dim}_\C {\Cal O}_{\L}(\Delta^n,p)\in \N$ is upper semicontinuous.
\endroster
\endproclaim

\noindent
As usual, we can now derive from Theorem~7.5 a complex algebraic
Frobenius theorem. As the dimension of orbits is an integer depending
upper-semi-continuously on the point $p$, it assumes its
maximal value $:=e$ at a Zariski-generic point $p\in\Delta^n$.
Let $\Cal V_M(p)$ be the intersection of a small polydisc 
$\Cal V_{\C^n}(p)$ with $M$.

\proclaim{Corollary~7.6} 
Let $p\in \Delta^n$ with $\text{\rm dim}_\C \, {\Cal
O}_{\L}(\Delta^n,p) =\max_{q\in {\Cal V}_{\Delta^n}(p)} \text{\rm
dim}_\C {\Cal O}_{\L}(\Delta^n,q)$ $=:e$. Then ${\Cal V}_{\Delta^n}(p)$ is
algebraically foliated by $\L$-integral manifolds of dimension $e$.
\endproclaim

\noindent
Using the properties summarized in~\thetag{7.3}, it is also easy to 
show that there exists a proper algebraic subvariety of $\Delta^n$
outside of which the dimension of $\L$-orbits is maximal equal to $e$.
Finally, we can derive from the previous considerations a real
algebraic CR-foliation theorem, useful to prove Theorem~2.7.

\proclaim{Corollary~7.7} 
Suppose that $p\in M$ is a point at which the dimension of nearby CR
orbits is maximal and locally constant equal to $2m+d-e$.  Then a
neighborhood ${\Cal V}_M(p)$ is real algebraically foliated by
CR-orbits and there exist $e$ holomorphic algebraic functions
$h_1,\ldots,h_e$ with $\partial h_1\wedge \cdots \wedge
\partial h_e(p)\neq 0$ such that
\roster
\item"{\bf (1)}"
$M$ is contained in $\{h_1=\bar h_1,\ldots,h_e=\bar h_e\}$.
In other words, $M$ is contained in a transverse intersection
of $e$ Levi-flat hypersurfaces in general position.
\item"{\bf (2)}"
Each manifold $M_c= M\cap \{h_1 = c_1,\ldots, h_e =
c_e\}$ is a CR-orbit of $M$.
\endroster
\endproclaim

\demo{Proof} 
As in the proof of Proposition~6.8 {\bf (2)}, we can intersect the
foliation provided by Corollary~7.6 with the antidiagonal
$\underline{\Lambda}$ to produce the algebraic foliation of $M$ by its
CR orbits. Let therefore $h_1,\ldots,h_e$ be real algebraic functions
over $(M,p)$ with linearly independent differential such that the
level-sets $\{h_1=c_1,\ldots, h_e=c_e\}$, $c_j\in\R$, are the plaques
of this foliation.  Since the $h_j$'s are constant on CR-orbits, they
are CR on $M$. Consequently, they extend as holomorphic algebraic
functions in a neighborhood of $(M,p)$ in $\C^n$, by the theorem of
Severi-Tomassini.
\qed
\enddemo

\noindent
Granted Corollary~7.7, a nowhere minimal CR-generic $M$ is contained
in at least one Levi-flat algebraic hypersurface in a neighborhood of
a Zariski-generic point. We are now in position to prove Theorem~2.7.

\proclaim{Theorem~7.8}
Let $M$ be a real algebraic CR-generic manifold in $\C^n$, let $p\in
M$, assume that $(M,p)$ is contained in the Levi-flat hypersurface
$\{z_n=\bar z_n\}$ and that $M$ possesses a nontrivial
infinitesimal CR-automorphism inducing an algebraic foliation.  Then
$M$ admits a local one parameter family of nonalgebraic biholomorphic
self-maps.
\endproclaim

\demo{Proof}
Let $\Cal X$ be a $(1,0)$-vector field with holomorphic algebraic
coefficients in a neighborhood of $p$ and let $\Cal K:=\hbox{Re}(\Cal
X)$ be the associated infinitesimal CR-auto\-morphism, which is tangent
to $(M,p)$. Since by assumption, the foliation induced by $\Cal X$ is
complex algebraic, we can assume after perharps multiplying $\Cal X$
by a nonzero algebraic function that its complex flow $(u,z)\mapsto
\exp (u\Cal X)(z)=:\varphi(u,z)$ is algebraic, {\it see} Lemma~9.12
below. It is well known that for $u$ real, the flow of $\Cal X$
coincides with the real flow of $\Cal K=\hbox{Re} \,
\Cal X$, hence it stabilizes $(M,p)$. 
By definition, $\partial_u \varphi(u,z)=\Cal X(\varphi(u,z))$. Since $\Cal
X\neq 0$, then obviously $\partial_u \varphi(u,z)\not\equiv 0$.  Let now
$\varpi(z_n)$ be an arbitrary nonzero holomorphic {\it nonalgebraic}
function with $\varpi(0)=0$ which is {\it real}, {\it i.e.}
$\overline{\varpi}(z_n)\equiv
\varpi(z_n)$ and satisfies $\partial_u\varphi(\varpi(z_n),z)\not\equiv 0$.
Then the map $z\mapsto \exp (\varpi(z_n) \Cal
X)(z)=\varphi(\varpi(z_n),z)$ is a holomorphic nonalgebraic
biholomorphism, because, if it where algebraic, the two conditions\,:
$\varphi(\varpi(z_n),z)$ is algebraic and
$\partial_u\varphi(\varpi(z_n),z)\not\equiv 0$ would clearly imply that
$\varpi$ is algebraic. Finally, this map sends $M$ into $M$, because
$\text{\rm Re} \, \Cal X$ is an infinitesimal CR-automorphism of $M$
{\it and because $\varpi$ is real on $M$}, which is the crucial point.
To get a one-parameter family of such maps, just take $\varphi(s \, 
\varpi(z_n), z)$ with $s\in\R$ small.
\qed
\enddemo

\head \S8. Partial algebraicity\endhead

\subhead 8.1.~Propagation of algebraicity \endsubhead 
Let $\L=\{L^{\alpha}\}_{\alpha\in A}$ be a set of vector fields
defined over $\Delta^n$ in $\C^n$ with complex algebraic flow.
Without loss of generality, we can assume that $k_\alpha=1$ for all
$\alpha\in A$. Following our scheme of proof for Theorem~2.5, it
remains now only to prove Theorem~5.2 {\bf (2)}, namely\,:

\proclaim{Theorem 8.2} 
A holomorphic function $g\in {\Cal H}(\Delta^n,\C)$ is algebraic on each
$\L$-orbit if and only it is algebraic on each $\L$-integral curve.
\endproclaim

\demo{Proof} 
We can assume $\L\neq \{0\}$. Of course, any small open piece of an
$\L$-orbit being algebraically locally equivalent to some
$\Delta^e$ by part {\bf (1)} of Theorem~5.2, 
part {\bf (2)} is then reduced to the case where $U$ is a
single $\L$-orbit, {\it i.e.} $U$ is $\L$-minimal.
\qed
\enddemo

\proclaim{Theorem 8.2'}
Assume that $\Delta^n$ is $\L$-minimal at $0$
and let $g\in {\Cal H}(\Delta^n,\C)$. If $g$ is algebraic
on each $\L$-integral curve, then $g$ is algebraic over $\Delta^n$.
\endproclaim

\remark{Remark} 
This theorem is proved by Sharipov-Sukhov in [14] in case there
exist $L^1,\ldots,L^n \in \L$ and $p\in U$ such that $\text{\rm rk}_\C \,
(L^1(p),\ldots,L^n(p))=n$, an assumption which corresponds essentially
to Segre-transversality of $M$ [5] in $\C^n$ 
in place of the more general assumption of
$\{\Cal L,\underline{\Cal L}\}$-orbit
minimality in the complexification ${\Cal M}$.
\endremark

\smallskip

We denote by $U$ a neighborhood of $0$ in $\Delta^n$ such that at
every point $q\in U$, the local $\L$-orbit of $q$ in $U$ contains an
open neighborhood of $q$. We shall say that $U$ is $\L$-minimal
(locally). Hence for each germ $\Lambda\subset U$ of a $\C$-algebraic
manifold with $\text{\rm dim} \
\Lambda < n$ and $\Lambda\neq \emptyset$, there exist $p\in \Lambda$
and $L\in \L$ such that $L(p)\not\in T_p\Lambda$. Otherwise, $\Lambda$
would be $\L$-integral, in contradiction with $\L$-minimality of $U$.
Of course, $\{q\in \Lambda \: L(q) \in T_q \Lambda \
\forall \ L \in \L\}$ is a proper closed $\C$-algebraic subset of
$\Lambda$. More generally, we shall need to consider this non-tangentiality
property along some families of such manifolds $\Lambda$ which 
foliate subdomains of $\Delta^n$.

\subhead 8.3.~Complex algebraic foliations \endsubhead
A regular holomorphic foliation ${\Cal F}$ of a subdomain $V\subset
\C^n$ is called {\it algebraic} if its transition maps are
$\C$-algebraic.  For short, we shall say in the sequel that ``${\Cal
F}$ is an ${\Cal A}$-foliation''.  We denote by $\Cal A(X,\C)$ the
ring of holomorphic algebraic functions on the complex manifold $X$.
Let now $\Lambda_{\Cal F}(p)$ be the leaf of ${\Cal F}$ through $p$
and let $m:=\text{\rm dim}_\C \, {\Cal F}$. We will consider only
local foliations, so that there is no restriction to assume that all
neighborhoods ${\Cal V}_{\C^n}(q)$, $V$, $W$, {\it etc.} in the sequel
are {\it foliation boxes}, which will simplify our
considerations. Recall that this means that ${\Cal F}$ is represented
in an open set, say $V$, by an algebraic (global) coordinate system
$\Phi\: V\to \Delta^m\times \Delta^{m-n}$, with respect to which the
leaves of $\Cal F$ in $\Delta^n$ are represented by the ``plaques''
$\Delta^m\times\{\zeta_{(n-m)}\}$. Thus, in $V$, the leaves of
$\Lambda_{\Cal F}(q)$ of the $\Cal A$-foliation $\Cal F$ are simply
the preimages $\Phi^{-1} (\Delta^m \times
\{\zeta_{(n-m)}(q)\})$, if we denote $\Phi(q)=\zeta_{(m)}(q)\times
\zeta_{(n-m)}(q)$. Let now $p\in V$, and choose
an arbitrary complex algebraic $(n-m)$-dimensional submanifold
$H\subset V$ passing through $p$ and satisfying $T_pH\oplus
T_p\Lambda_{\Cal F} (p)=T_pV$. Thus $\text{\rm dim}_\C \, H+\text{\rm
dim}_\C \, {\Cal F}=n$. Then there exists a neighborhood $W:={\Cal
V}_{\C^n}(p)$ such that $\cup_{q\in H}
\Lambda_{\Cal F} (q)\cap W={\Cal F}|_W$. Again, $W$ is a foliation
box here and the leaves $\Lambda_{\Cal F}(q)$ are embedded closed
$\C$-algebraic submanifolds of $W$ which are the plaques of $\Cal
F\v_W$. Finally, recall that to any vector field $L\in \L$ with $L\neq
0$ and to any point $p\in U$ with $L(p)\neq 0$, there is associated a
neighborhood $V={\Cal V}_{\C^n}(p)$ which is algebraically foliated by
$L$-integral curves. Such foliations will be denoted by $\Cal F_L$ in
the sequel.

\subhead 8.4.~Description of the proof of Theorem~8.2'\endsubhead
We now introduce an important notation. We shall
write $g\in {\Cal A}_{\Cal F}(V)$ if $g$ restricted to each
leaf of ${\Cal F}|_V$ is algebraic. To establish Theorem~8.2', we will prove
the following statement by induction on the integer $m\in\N$\,:

\smallskip
\roster
\item"$(*)$" 
For all $m$, $1\leq m \leq n$, there exist $p\in U$ and an
$m$-dimensional foliation ${\Cal F}$ of $V={\Cal V}_{\C^n}(p)$ such that
$g\in {\Cal A}_{\Cal F}(V)$.
\endroster

\noindent
Then the desired Theorem~8.2' will be just $(*)$ for $m=n$. We
already know $(*)$ for $m=1$, since $g\in {\Cal A}_{\Cal F}({\Cal
F}_L)$ for every nonzero vector field $L\in
\L$. Let $m\leq n-1$. We have to assume $(*)$ for $m$, {\it i.e.} $(*)_m$ and
to deduce $(*)$ for $m+1$, {\it i.e.} $(*)_{m+1}$.
{\it We can describe now this inductional implication in the large, 
except some two main lemmas that will be rejected below.}
First, by $\L$-minimality of $U$, we have obviously\,:

\proclaim{Lemma~8.5}
Let $p\in U$, $V={\Cal V}_{\C^n}(p)$, and let ${\Cal F}$ be an ${\Cal
A}$-foliation of $V$ with $\text{\rm dim}_\C \, {\Cal F}=m$, $1\leq
m\leq n-1$. Then for all $q\in V$, there exist $r\in
\Lambda_{\Cal F}(q)$ arbitrarily close to
$q$ and $L\in \L$ such that the vector $L(r)$ is not tangent 
to the leaf $\Lambda_{\Cal F}(q)$.
\endproclaim

Indeed, otherwise, the manifold $\Lambda_{\Cal F}(q)$ would be an
$\L$-integral manifold in a neighborhood of $q$ of positive
codimension, a contradiction with local $\L$-minimality of $U$ at
every point. Let now $z=(z_1,\ldots,z_n)$ denote algebraic coordinates
on $U$. Here is our first main lemma, inspired from [14], which will
be the first step of our proof of the inductional implication.  This
lemma will state that if $g$ is algebraic on the leaves of ${\Cal
F}|_V$, then $g$ and all its derivatives are algebraic on the leaves
of ${\Cal F}|_W$ after restriction to a possibly smaller subdomain
$W\subset V$.

\proclaim{Lemma~8.6} 
If $g\in {\Cal A}_{\Cal F}(V)$, then there exists a nonempty subdomain
$W\subset V$ such that all partial derivatives
$\partial^{\beta}_z g\in {\Cal A}_{\Cal F}(W)$ too,
for all $\beta\in \N^n$.
\endproclaim

It is easy to 
check that this condition is independent of coordinates.

\proclaim{Lemma~8.7} 
Let $V\subset \C^n$ be a domain, let $z$ be complex algebraic
coordinates on $V$,
let ${\Cal F}$ be an ${\Cal A}$-foliation on $V$, let $g\in {\Cal
H}(V,\C)$ and let $\Phi\: V\to
\Phi(V)$ be an ${\Cal A}$-biholomorphism, $w=\Phi(z)$. If
$\partial_z^{\beta}g\in {\Cal A}_{\Cal F}(V)$ $\forall \ \beta \in
\N^n$, then $\partial^{\beta}_w(g\circ \Phi^{-1}) \in {\Cal
A}_{\Phi_*{\Cal F}}(\Phi(V))$, $\forall \ \beta \in \N^n$.
\endproclaim

\demo{Proof} 
Simple application of the chain rule, because
there exist universal polynomials $P_{\beta}$ with
$\partial_w^{\beta}(g\circ
\Phi^{-1})=P_{\beta}(\{\partial_w^{\gamma}\Phi^{-1}\}_{\gamma \leq
\beta}, \{(\partial_z^{\gamma}g)\circ \Phi^{-1}\}_{\gamma\leq
\beta})$.
\qed
\enddemo

To achieve the second step of the proof, we have to construct a
foliation ${\Cal F}^{1+}$ with $\text{\rm dim} {\Cal F}^{1+}=m+1$
satisfying $(*)_{m+1}$. Let $p$, $V$ and ${\Cal F}$ be as in $(*)_m$
and let $q\in V$ be arbitrary {\it with $V$ like the $W$ in
Lemma~8.6} (shrinking $V$ if necessary).  Choose now $r\in
\Lambda_{\Cal F}(q)$ as in Lemma~8.5 with $L(r)\not\in
T_r\Lambda_{\Cal F}(q)$. Changing notation, we will now denote this
$r$ by $p$. Further, let $H$ be a piece of an ${\Cal A}$-manifold
through $p$ with $T_pH \oplus T_p \Lambda_{\Cal F}(p)=T_pV$, so
$\text{\rm dim}_\C \, H=n-m$. Then $\Cal F\v_V=\cup_{q\in H}
\Lambda_{\Cal F}(q)\cap V$ (after shrinking $V$ and still with the
convention about small open sets being foliation boxes), where, of
course, $\Lambda_{\Cal F}(q) \cap \Lambda_{\Cal F}(q') =\emptyset$ if
$q\neq q'$. Now, instead of $H$, let us choose an arbitrary piece
$H^{1+}$ of an ${\Cal A}$-manifold with $p\in H^{1+}$ and
$T_pH^{1+}\oplus \C L(p)\oplus T_p \Lambda_{\Cal F}(p)= T_pV$, so
$\text{\rm dim}_\C H^{1+}=n-m-1$. To the triple $({\Cal F}, L,
H^{1+})$, we can finally associate an ${\Cal A}$-foliation ${\Cal
F}^{1+}= {\Cal F}^{1+}({\Cal F}, L, H^{1+})$ of $\Cal V_{\C^n}(p)$
with $\text{\rm dim}_\C \, {\Cal F}^{1+}=m+1$, which is defined as
follows and will be the important object to get $(*)_{m+1}$. Simply,
this foliation will be constructed by flowing the leaves of $\Cal F$ a
bit along the nontangential flow lines of $L$, thus gaining one unit
in dimension. Precisely\,:
\roster
\item"{\bf 1.}" 
The leaves $\Lambda_{{\Cal
F}^{1+}}(q):=\{L_s(r) \: \v s \v < \delta, r\in \Lambda_{\Cal
F}(q)\}$, for $q\in H^{1+}$, 
\item"{\bf 2.}" 
${\Cal F}^{1+}=\cup_{q\in H^{1+}}
\Lambda_{{\Cal F}^{1+}}(q)$. 
\endroster
By construction, 
the foliation ${\Cal F}^{1+}$ clearly has the following structure\,:
\roster
\item"{\bf 3.}" 
${\Cal F}_L$ is a subfoliation of ${\Cal F}^{1+}$ in ${\Cal V}_{\C^n}(p)$.
\item"{\bf 4.}" 
${\Cal F}$ is not in general a subfoliation of ${\Cal F}^{1+}$, but
each leaf of ${\Cal F}^{1+}$ contains at least one leaf of ${\Cal F}$\,:
$\Lambda_{{\Cal F}^{1+}}(q) \supset \Lambda_{\Cal F}(q)$ for $q\in
H^{1+}$.
\endroster

\remark{Important remark} 
Such a step-by-step construction of several foliations is necessary
because, even if there may exist $n$ linearly independent vector
fields $L^1, \ldots, L^n$ in $\L$ over $U$, there might not exist a
coordinate system whose coordinate line correspond to the flow lines
of the $L^j$'s. More specifically, and for the same reason, even in
the above construction of $\Cal F^{1+}$, there does not exist in
general a system of coordinates $(x_1,\ldots,x_m,
x_{m+1},y_1,\ldots,y_{n-m-1})$ such that the $x_{m+1}$-lines
correspond to the integral curves of $L$ and such that,
simultaneously, the $(x_1,\ldots,x_m)$-planes correspond to leaves of
${\Cal F}$. In particular, one cannot in fact reduce our Theorem~8.2'
directly to a classical theorem due to Bochner-Martin (1949) [3] which
states that given $g\in {\Cal H}(\Delta^n, \C)$ algebraic on
$z_i$-lines, $i=1,\ldots,n$, then $g$ is algebraic. To argue this
general impossibility, let us consider for instance for $n=2$, some
two linearly independent vector fields $L_1:=\partial / \partial x$
and $L_2:=\partial /
\partial y+a(x,y) \partial / \partial x$ (general form of such a pair
of vector field after a convenient choice of coordinates $(x,y)$).  In
this case, it is easy to see that there further exists a
biholomorphism of ${\Cal V}_{\C^2}(0)$ transforming $L_1$ and $L_2$
into $\partial / \partial x$ and $\partial / \partial y$. However,
this fails if $n=3$.  Indeed, let $L_1=\partial / \partial x$ and
$L_2=\partial / \partial y + x\, \partial / \partial z$ in
$\C^3(x,y,z)$. Then there does not exist $\Phi=(\Phi_1,\Phi_2,\Phi_3)$
a biholomorphism of ${\Cal V}_{\C^3}(0)$ transforming $L_1$ in
$\partial /\partial x$ and $L_2$ in $\partial /\partial y$, because
such $\Phi$ would have to satisfy $\Phi_{3x}=0$ and
$\Phi_{3y}+x\Phi_{3z}=0$, whence $\Phi_{3z}=0$, $\Phi_{3y}=0$,
$\Phi_3=ct.$, a contradiction.
\endremark

\smallskip
The lemma below is analogous to Lemma~8.7 and is also elementary.

\proclaim{Lemma~8.8} 
Let $g\in {\Cal H}(V)$, $\partial_z^{\beta}g\in {\Cal A}_{\Cal F}
(V)$, $\forall \ \beta \in
\N^n$. Let $q\in H^{1+}$, let $\Lambda:=\Lambda_{\Cal F}(q)$, let
$\Lambda^{1+}= \Lambda_{{\Cal F}^{1+}}(q)$ $(\supset \Lambda)$. Let
$z^+$ be ${\Cal A}$-coordinates on $\Lambda^{1+}$ and put
$g^+=g|_{\Lambda^{1+}}$. Then $\partial_{z^+}^{\beta}(g^+|_{\Lambda})$
is algebraic on $\Lambda$, for all $\beta\in \N^{m+1}$.\qed
\endproclaim

Our second main lemma will be as follows. To understand it concretely,
the reader may read parallely its formulation in coordinates given in
Lemma~8.12 below. With our notations, in an arbitrary fixed leaf
$\Lambda_{\Cal F^{1+}}(q)$, this lemma will state that the restriction
of $g$ to $\Lambda_{\Cal F^{1+}}(q)$ is algebraic provided it
is algebraic on every ``vertical $1$-dimensional leaf'' $\Lambda_{\Cal
F_L}(r)\subset\Lambda_{\Cal F^{1+}}(q)$, $r\in \Lambda_{\Cal F}(q)$
and provided all its jets of any order are algebraic, when restricted
to the fixed ``horizontal'' $m$-dimensional manifold $\Lambda_{\Cal
F}(q)\subset \Lambda_{\Cal F^{1+}}(q)$.

\proclaim{Lemma~8.9} 
Let $p$, $V$, ${\Cal F}$, ${\Cal F}^{1+}$ and $L$ be as above and let
$z_q^+\in \C^{m+1}$ denote some algebraic local coordinates on $\Lambda_{{\Cal
F}^{1+}}(q)$. Let $g\in {\Cal H}(V,\C)$. Then $g\in {\Cal A}_{{\Cal
F}^{1+}}(V)$, {\it i.e.} its restriction $g|_{\Lambda_{\Cal
F}^{1+}(q)}$ is algebraic for all $q\in H^{1+}$, provided it
satisfies the following two conditions for each $q$\,:
\roster
\item"{\bf 1.}"
Algebraicity on the flow lines of $L$\,:
$g\in {\Cal A}_{{\Cal F}_L}(V)$.
\item"{\bf 2.}"
Algebraicity of all the differentials of $g$ after restriction 
to the central one-codimensio\-nal leaf $\Lambda_{\Cal F}(q)$\,:
$\partial_{z_q^+}^{\beta} (g|_{\Lambda_{\Cal F}^{1+}(q)})
|_{\Lambda_{\Cal F}(q)}\in {\Cal A}(\Lambda_{\Cal F}(q),\C)$,
$\forall \ \beta\in \N^{m+1}$.
\endroster
\endproclaim

\remark{Remarks} 
{\bf 1.}~Notice that the first main
assumption of this lemma is already in the hypotheses of Theorem~5.2,
and that the second one follows from Lemma~8.6.

{\bf 2.}~In coordinates $z_q^+=(z_q,w_q)$ with $\Lambda_{\Cal F}(q)=
\{w_q=0\}$, it suffices in fact to require above that all the
$\partial^k_{w_q} g(z_q, 0)$ are algebraic, since then
$\partial^{\beta}_{z_q} \partial^k_{w_q} g(z_q,0)$ are algebraic, by
the stability of algebraicity under differentiation.  In other words,
it suffices to require only that the transversal jets are algebraic,
since the horizontal jets are then automatically algebraic.
\endremark

\smallskip

In conclusion, the foliation ${\Cal F}^{1+}$ is the sought foliation
satisfying $(*)_{m+1}$. Theorem~8.2' will then be proved once we have
proved Lemmas~8.6 and 8.9.
\qed 

\subhead 8.10.~Proofs of Lemmas~8.6 and 8.9\endsubhead 
It is now possible to reformulate our two main lemmas. Passing to
coordinates, there is given an open set ${\Cal U}={\Cal S}\times {\Cal
T} \subset
\C^a \times \C^b$, with, say ${\Cal S}=\Delta^a$, ${\Cal T}=\Delta^b$,
with $a\in \N_*$, $b\in \N_*$, $0\in {\Cal S}$, $0\in {\Cal T}$,
$(s,t)\in {\Cal S}\times {\Cal T}$, and a holomorphic function $g\:
(s,t)\mapsto g(s,t)$ defined over $\Cal S\times
\Cal T$. Here, we shall choose
$a:=m$, $b:=n-m$ for Lemma~8.6 and $a:=m$, $b:=1$ for Lemma~8.9.
Let $c:=a+b$.

\proclaim{Lemma~8.11} 
If $[{\Cal S} \ni s\mapsto g(s,t)\in \C] \in {\Cal A}({\Cal S}, \C)$,
$\forall \ t\in {\Cal T}$, then there exists a nonempty open subset
${\Cal T}_1\subset {\Cal T}$ such that $[{\Cal S} \ni s\mapsto
\partial^{\gamma}_{s,t} g(s,t)\in \C ]\in {\Cal A}({\Cal S}, \C)$,
$\forall \ t\in {\Cal T}_1$, $\forall \ \gamma\in \N^c$.
\endproclaim

\remark{Remark} 
For Lemma~8.6, the sets ${\Cal S}\times \{t\}$ play the r\^ole 
of the leaves of $\Cal F$.
\endremark

\proclaim{Lemma~8.12} 
If $[{\Cal S}\ni s\mapsto
\partial^{\beta}_tg(s,0)\in \C ]\in {\Cal A}({\Cal S},\C)$ $\forall \
\beta \in \N^b$ and if $[{\Cal T} \ni t\mapsto g(s,t)\in \C ]\in {\Cal
A}({\Cal T}, \C)$, $\forall \, s\in {\Cal S}$, then $g(s,t)\in {\Cal
A}({\Cal S}\times {\Cal T}, \C)$.
\endproclaim

\remark{Remark} 
For Lemma~8.9, after fixing $q\in H^{1+}$, the set ${\Cal S} \times
\{0\}$ is $\Lambda_{\Cal F} (q)$, the sets $\{s\}\times {\Cal T}$ are
leaves of ${\Cal F}_L$ contained in $\Lambda_{\Cal F}^{1+}(q)$ and
${\Cal S}\times {\Cal T}$ is $\Lambda_{\Cal F}^{1+}(q)$.
\endremark

\demo{Proof of Lemma~8.11}
By assumption, $\forall \, t\in {\Cal T}$, $\exists \, N_t\in \N$, $N_t
\geq 1$, $\exists \, A_t\in \N$, $A_t \geq 1$, $\exists \, a_{i\alpha,t}\in
\C$, $1\leq i \leq N_t$, $|\alpha|\leq A_t$, $\alpha \in \N^a$ and
{\it irreducible} polynomials
$$ 
\C[X,s]\ni
P_t(X,s):=
\sum_{i=0}^{N_t} \sum_{|\alpha|\leq A_t} a_{i\alpha,t} \ X^i \
s^{\alpha},
\tag 8.13
$$
such that $P_t(g(s,t),s)\equiv_s 0$.
Here, the notation $\Phi(s,t)\equiv_s 0$ means that the
formal power series $\Phi(s,t)$ vanishes identically when 
considered as a series in $s$ only. Let us transform
the $a_{i\alpha,t}$ first to make them depend in a nice 
way with respect to $t$.

\proclaim{Lemma~8.14} 
Then there exists a nonempty open set ${\Cal T}_1 \subset \subset
{\Cal T}$ and an irreducible polynomial of uniformly bounded degree
having coefficients $a_{i\alpha}(t)$ holomorphic over $\Cal T_1$, namely
there exists\,:
$$ 
{\Cal H}({\Cal T}_1)[X,s]\ni
P(X,s;t):=\sum_{i=0}^N\sum_{|\alpha|\leq
A} a_{i\alpha}(t) \ X^i \ s^{\alpha}, \ \ \ \ \ a_{i\alpha}\in {\Cal
H}({\Cal T}_1),
\tag 8.15
$$ 
such that $P(g(s,t),s; t)\equiv_{s,t} 0$, for all $s\in {\Cal S}$
and $t\in {\Cal
T}_1$.
\endproclaim

\demo{End of proof of Lemma~8.11}
Taking Lemma~8.14 for granted, it then suffices to differentiate
$P(g(s,t),s;t)\equiv 0$ with respect to $t$, namely to apply
$\partial_t^\beta$, $\v \beta\v=1$, and to eliminate $g(s,t)$ from the
system of algebraic equations\,:
$$
\left\{
\aligned
&
\sum_{i=0}^N \sum_{|\alpha|\leq A}\left( \partial_t^\beta 
a_{i\alpha}(t)\right. g(s,t)+
\left.
i \, a_{i\alpha}(t) \, \partial_t^\beta g(s,t)\right) \, g(s,t)^{i-1} 
\, s^{\alpha}=0,\\
& 
P(g(s,t),s; t)=0
\endaligned\right.
\tag 8.16
$$ 
this giving the algebraicity of $s\mapsto \partial_t^\beta g(s,t)$ for
all $\v \beta\v=1$. General induction is analogous\,: simply replace
$g(s,t)$ by $\partial_t^\beta g(s,t)$.
\qed 
\enddemo
\enddemo

\demo{Proof of Lemma~8.14}
To prove the assertion, consider the countable many sets ${\Cal
E}_{N,A}=\{t\in {\Cal T}; $ $N_t= N, A_t= A\}$ whose union
$\bigcup_{N\geq 1, A\geq 1} {\Cal E}_{N,A}$ equals ${\Cal T}$, 
by~\thetag{8.13}. Of course,
the union of their closures $\bigcup_{N\geq 1, A\geq 1} \overline{{\Cal
E}_{N,A}}$ equals ${\Cal T}$ as well. Thanks to Baire's category 
theorem, at least one closure
$\overline{{\Cal E}_{N,A}}$ has nonempty interior. 
This is the main trick. Thus, there exists a nonempty open
polydisc $\Cal T_1\subset \subset {\Cal T}$ with center $t_1\in\Cal T_1$
such that ${\Cal T}_1 \subset \overline{{\Cal E}_{N,A}}$. Then
we have polynomial relations\,:
$$
P_t(g(s,t),s)=\sum_{i=0}^N\sum_{|\alpha|\leq A} a_{i\alpha,t} \ g(s,t)^i \
s^{\alpha} \equiv_s 0, \ \ \ \ \ \text{\rm for}
\ \ t\in {\Cal T}_1 \cap {\Cal E}_{N,A},
\tag 8.17
$$ 
with a uniform bound for the degrees on a dense subset of
$\Cal T_1$. This means that the finite set of
functions
$$ 
{\Cal C}_t=\{{\Cal S}\ni s\mapsto g(s,t)^i \, s^{\alpha}\in \C \}_{0\leq
i\leq N, |\alpha|\leq A}
\tag 8.18
$$ 
is always linearly dependent for $t\in {\Cal T}_1 \cap {\Cal
E}_{N,A}$, and in particular for $t=t_1$. Let us 
normalize one $a_{i_*\alpha_*,
t_1}=1$, for some $i_*, \alpha_*$ and let us denote, after renumbering its
elements, the set $\Cal C_t$ by ${\Cal
C}_t:=\{{\Cal S} \ni s\mapsto h_k(s,t)\in \C
\}_{1\leq k \leq B}$, with $B\in \N_*$ and $h_B(s,t)=:g^{i_*}(s,t)
\, s^{\alpha_*}$. 
Further, let us develope each $h_k$ in power series with respect
to $s$\,:
$$ 
h_k(s,t)=\sum_{l\in \N^a} h_{kl}(t) \, s^l, \ \ \ \ \ h_{kl}\in {\Cal
H}({\Cal T}_1).
\tag 8.19
$$ 
\proclaim{Lemma~8.20}
Now, the following properties hold\,:
\roster
\item"{\bf (1)}"
The functions $h_1(s,t_1),\ldots,h_{B-1}(s,t_1)$ are linearly 
independent over $\C$.
\item"{\bf (2)}"
There exist $(B-1)$ pairwise distinct multiindices 
$l_{*1},\ldots,l_{*B-1}\in\N^a$ such that the determinant 
$\text{\rm det} \, ((h_{kl_{*j}}(t_1))_{1\leq k,j\leq B-1})$ is nonzero.
\item"{\bf (3)}"
Consequently, the functions of $s$, 
$h_1(s,t),\ldots,h_{B-1}(s,t)$ are linearly 
independent over $\C$ for all $t$ running in a neighborhood of $t_1$.
\endroster
\endproclaim

\demo{Proof}
Suppose by contradiction that there exist $d_1,\ldots,d_{B-1}\in\C$
not all zero such that $d_1h_1(s,t_1)+\cdots+d_{B-1}h_{B-1}(s,t_1)
\equiv_s 0$. Reexpressing the $h_k$'s in terms of $g$ and $s$, 
this means that there exists a nonzero polynomial $Q(X,s)$ such that
$Q(g(s,t_1),s)\equiv_s 0$. As $P_{t_1}(X,s)$ was assumed to be
irreducible, it follows that there exists a nonzero polynomial
$u(X,s)$ such that $Q(X,s)\equiv u(X,s)P_{t_1}(X,s)$. But the
monomials of $Q$ are exactly the same as those of $P_{t_1}$, except
one missing term $X^{i_*} \, s^{\alpha_*}$. For reasons of degree, $u$
is then a nonzero constant. This contradicts the fact that $P_{t_1}$
incorporates the monomial $X^{i_*} \, s^{\alpha_*}$, since 
$a_{i_*,\alpha_*,t_1}=1$ by our previous choice.
Thus, {\bf (1)} is proved. Then the $(B-1)\times\infty$-matrix of complex
coefficients $(h_{kl}(t_1))_{1\leq k\leq B-1, \, l\in\N^a}$ possesses
a nonzero $(B-1)\times (B-1)$ minor, which yields {\bf (2)} and then
{\bf (3)} evidently.
\qed
\enddemo

\demo{End of proof of Lemma~8.14}
Shrinking $\Cal T_1$ if necessary, we can therefore assume that the
functions $h_1(s,t),\ldots,h_{B-1}(s,t)$ are linearly independent for all
$t\in\Cal T_1$. Remembering that $h_1(s,t),\ldots,h_B(s,t)$ are, on the
contrary, linearly {\it dependent} for all $t\in\Cal T_1\cap \Cal E_{N,A}$,
we obtain that there exist $a_{1,t},\ldots,a_{B-1,t}\in\C$ such that
$$
a_{1,t} \, h_1(s,t)+\cdots+a_{B-1,t} \, h_{B-1}(s,t)\equiv_s h_B(s,t), 
\ \ \ \forall \, t\in\Cal T_1\cap \Cal E_{N,A},
\tag 8.21
$$
or equivalently, using the developement~\thetag{8.19},
$$
a_{1,t} \, h_{1\, l}(t)+\cdots+a_{B-1,t} \, h_{B-1 \, l}(t)=h_{B \, l}(t), 
\ \ \ \forall \, t\in\Cal T_1\cap \Cal E_{N,A}, \ \forall \, l\in\N^a.
\tag 8.22
$$
Writing in particular these equations for $l=l_{1*},\ldots,l_{B-1 \,
*}$, using property {\bf (2)} of Lemma 8.20 and Cramer's rule, we
deduce that the coefficients $a_{1,t},\ldots,a_{B-1, t}$ can be uniquely
expressed as rational functions with respect to the
$h_{k,l_{*j}}(t)$, $k=1,\ldots,B$, $j=1,\ldots,B-1$, with the
denominator $\text{\rm det} \, ((h_{kl_{*j}}(t))_{1\leq
k,j\leq B-1})$, nonvanishing over $\Cal T_1$. 
We thus have got holomorphic functions $a_1(t),\ldots,
a_{B-1}(t)$ over $\Cal T_1$ which are rational 
in the coefficients $h_{k,l_{*j}}(t)$ and which satisfy
$$
a_1(t) \, h_1(s,t)+\cdots+a_{B-1}(t) 
\, h_{B-1}(s,t)\equiv_s \, h_B (s,t),
\tag 8.23
$$
for all $t\in\Cal T_1\cap \Cal E_{N,A}$ and then for all $t\in \Cal T_1$,
by continuity. Reexpressing the $h_k$'s in terms of 
$g$ and $s$, we finally get~\thetag{8.15}.
\qed
\enddemo
\enddemo

\demo{Proof of Lemma~8.12} 
Direct differentiations with respect to $t$ yield the following more
explicit version of~\thetag{8.19}, where the r\^oles of $s$ and $t$ 
are exchanged\,:

\proclaim{Lemma~8.24}
There exist integers $c_{k,\beta,\gamma}$
and universal polynomials $\phi_{k,\beta,\gamma}$
such that each monomial $g(s,t)^k \, t^\beta$
can be expressed as
$$
g(s,t)^k \, t^\beta =\sum_{\gamma\in\N^b} \phi_{k,\beta,\gamma}
(\{\partial_t^\beta g(s,0) \}_{\v\beta\v\leq c_{k,\beta,\gamma}}) \, 
t^\gamma.
\tag 8.25
$$
\endproclaim

\demo{Proof}
Applying the differentiations $\partial_t^\gamma\v_{t=0}$
to $g(s,t)^k \, t^\beta$, we simply get 
$$
\phi_{k,\beta,\gamma} 
(\{\partial_t^\beta g(s,0) \}_{\v\beta\v\leq c_{k,\beta,\gamma}}):=
\partial_t^\gamma\v_{t=0}(g(s,t)^k \, t^\beta)/\gamma !.
\qed
\tag 8.26
$$
\enddemo

\noindent
This expression is very appropriate, because the partial derivatives
$\partial_t^\beta g(s,0)$ appear in an algebraic way. On the other
hand, the maps $t\mapsto g(s,t)$ are algebraic by assumption.
Inversing the r\^oles of $s$ and $t$, Lemma~8.14 then yields an open
subset $\Cal S_1\subset\subset \Cal S$ and an irreducible polynomial
$R(Y,t;s)=\sum_{k=0}^{N} \sum_{\beta\in\N^b, \v \beta\v \leq
B}b_{k\beta}(s) \, Y^k \, t^\beta$ such that $R(g(s,t),t;s)\equiv_{t,s}
0$. Here, the functions $b_{k\beta}(s)$, which are holomorphic over
$\Cal S_1$, are in fact {\it rational} with respect to a finite number
of the coefficients $\phi_{l,\beta,\gamma}(\{\partial_t^\beta
g(s,0)\}_{\v \beta\v
\leq c_{l,\beta,\gamma}})$ appearing in~\thetag{8.25}, 
as we have observed just after~\thetag{8.22} in our
application of Cramer's rule.  In summary, there exist an
integer $e\in\N_*$ and a polynomial relation
$$ 
R(g(s,t),t)=\sum_{\beta\in \N^b, \v\beta\v\leq B} \ \sum_{k=0}^{N}
\psi_{k\beta}(\{\partial_t^{\beta}g(s,0)\}_{\v\beta\v\leq e}) \ g(s,t)^k \,
t^\beta\equiv_{s,t} \, 0,
\tag 8.27
$$ 
with the $\psi_{k\beta}$ being rational. As by assumption, the partial
derivatives $\partial_t^{\beta}g(s,0)$ are all algebraic, using
elimination theory, we can transform~\thetag{8.27} into a polynomial
relation $S(g(s,t),s,t)=0$, where $S(Z,s,t)\in\C[Z,s,t]$, which shows
that $g(s,t)$ is algebraic. The proofs of Lemmas~8.14 and 8.12 are
complete.
\qed
\enddemo

In conclusion, Theorem~5.2 and Theorem~2.5 are now fully established. 
It remains now to explain what is the degree of algebraic degeneracy. \qed

\head \S9. Algebraic degeneracy \endhead

\subhead 9.1.~Real algebraic sets\endsubhead
We denote by $\C[z,\bar z]_\R\subset \C[z,\bar z]$ the ideal of {\it
real} polynomials $P(z,\bar z)$, namely those satisfying $P(z,\bar
z)\equiv \bar P(\bar z,z)$. Then $\C[z,\bar z]_\R$ is isomorphic to
$\R[\hbox{Re} \, z,\hbox{Im} \, z]$. Let $\Sigma\subset \C^n$ be a
real algebraic set, defined as the zero set of a collection of
elements of $\C[z,\bar z]_\R$. Let $\Cal J(\Sigma)$ denote the ideal
of polynomials vanishing on $\Sigma$. This ideal is prime if and only
if $\Sigma$ is irreducible, which we will suppose throughout \S9. The
field of fraction $K(\Sigma)$ of the entire ring $\C[z,\bar z]_\R/\Cal
J(\Sigma)$ is called the {\it field of rational functions on
$\Sigma$}. Its transcendence degree over $\C$ is called the {\it
dimension} of $\Sigma$. Let $\delta$ be this dimension. If
$P_1,\ldots,P_\sigma\in\C[z,\bar z]_\R$ is a system of generators of
the prime ideal $\Cal J(\Sigma)$, the generic rank of the complex
$\sigma\times 2n$ Jacobian matrix
$$
J_P(z,\bar z):= \left( {\partial P_i\over \partial z_k} (z,\bar z) \ \
{\partial P_i\over \partial \bar z_k} (z,\bar z) \right)_{ 
1\leq i\leq \sigma, \,1\leq k\leq
n},
\tag 9.2
$$
is then equal to $d:=n-\delta$ over $\Sigma$, the {\it codimension}
of $\Sigma$. 
The set $\Sigma_{reg}$ at which the
rank of $J_P(z,\bar z)$ is equal to $d$ does not depend on the
choice of a system of generators for $\Cal J(\Sigma)$ and is called
the {\it set of regular points of $\Sigma$}, in the algebraic
sense. Its complement $\Sigma_{sing}:=\Sigma\backslash \Sigma_{reg}$
is a proper real algebraic subvariety of $\Sigma$ called its {\it set
of singular points}.

\example{Example~9.3}
The celebrated {\it Whitney umbrella} is the cubic of $\R^3$ defined
by $W:=\{(x_1,x_2,x_3)\in\R^3\: x_3x_1^2= x_2^2\}$. In the algebraic
sense, $W_{sing}=0\times 0\times \R$ and the set $W_{reg}=W\cap
\{(x_1,x_2)\neq (0,0)\}\subset \{x_1\neq 0\}$ is not connected and not
dense in $W$ for the euclidean topology.  This example also shows that
the set of geometrically singular points of $W$, equal to $\{(0,0,x_3)
\: x_3\geq 0\}$ is only semi-algebraic.  Let us consider the tube
$\Sigma:= W\times (i\R)^3$ over $W$ in $\C^3$, namely
$\Sigma:=\{(z_1,z_2,z_3)\in\C^3\: x_3x_1^2=x_2^2\}$, $x_j=\hbox{Re} \,
z_j$. The set $\Sigma_{sing}=\{x_1=x_2=0\}= \R\times \R\times \C$ is
then $1$-algebraically degenerate and nowhere minimal. The set
$\Sigma_{reg}=\Sigma\cap \{x_1\neq 0\}$ is a rigid real algebraic
hypersurface globally defined as the graph $x_3=x_2^2/ x_1^2$ on
$\C^3\backslash \{x_1=0\}$ which carries the two CR vector fields
$\bar L_1:= x_1 \partial / \partial \bar z_1- 2x_3 \partial / \partial
\bar z_3$ and $\bar L_2:= x_1^2 \partial / \partial \bar z_2+
2x_2 \partial /\partial \bar z_3$. Since the Lie bracket $[L_2,\bar
L_2]=x_1^2\partial / \partial \bar z_3-x_1^2\partial /\partial z_3$ is
linearly independent with $(\bar L_1, \bar L_2)$ at every point, the
hypersurface $\Sigma_{reg}$ is minimal at every point.  Further, if $\Cal
T:=\sum_{j=1}^3 a_j(z)\partial /\partial z_j$ is a local holomorphic
algebraic vector field tangent to $\Sigma_{reg}$, we have
$a_1(z)x_1x_3-a_2(z)x_2+a_3(z)x_1^2/2=0$, whence
$a_1=a_2=a_3=0$, namely $\Sigma_{reg}$ is algebraically nondegenerate
(it is even Levi-nondegenerate at every point).
This contrasts with $\Sigma_{sing}$. Let us finally illustrate
Theorems~2.5,~2.6 and~3.1. Let $f\: (\Sigma,p)\to (\Sigma,p)$ be a
germ of biholomorphism fixing $p$. If $p\in \Sigma_{sing}\cap
\{x_3< 0\}$, there exists a local perturbation $\phi\circ f$ of
$f$ such that $\nabla^{tr}(\phi\circ f)=1$. If $p\in \Sigma\cap
\{x_3\geq 0\}$, then $f$ is necessarily algebraic.
\endexample

\subhead 9.4.~Extrinsic complexification \endsubhead We define
$\Sigma^c:=\{(z,\zeta)\in\C^{2n}\: P(z,\zeta)=0, \ \forall \, P\in\Cal
J(\Sigma)\}$. Then $\Sigma^c$ is a complex algebraic subset of
$\C^{2n}$ with $\Sigma^c\cap \underline{\Lambda} =\Sigma$, where 
$\underline{\Lambda}=\{\zeta=\bar z\}$, if we identify
$\{(z,\bar z)\: z\in \Sigma\}$ with $\{z\in \Sigma\}$. This set $\Sigma^c$
is useful because $\Sigma_{reg}$ need not be connected. Using the
rank property of the Jacobian~\thetag{9.2} and the connectedness of
the regular part of complex algebraic sets, one can check that
$\Sigma^c$ is irreducible if $\Sigma$ is. Further, using the
connectedness of the regular part $\Sigma_{reg}^c$ of $\Sigma^c$, one
can show that if there exists a point $p\in \Sigma$ such that
$r(p,\bar p)\neq 0$, where $r(z,\bar z)\in\C[z,\bar z]_\R$, then the
set $\{z\in \Sigma\: r(z,\bar z)\neq 0\}$ is dense in $\Sigma_{reg}$
for the euclidean topology. An arbitrary union of sets of the form
$\{z\in\Sigma\: r(z,\bar z)\neq 0\}$, where $r$ does not vanish
identically on $\Sigma$, will be called a Zariski open subset of
$\Sigma$. A point $p\in\Sigma$ will be called {\it Zariski-generic} if
it runs in some Zariski open subset of $\Sigma$. 

\subhead 9.5.~Intrinsic complexification \endsubhead
By noetherianity of $\C[z]$, the intersection of all complex algebraic 
subvarieties containing $\Sigma$ is a complex algebraic subvariety, called
the {\it intrinsic complexification} of $\Sigma$ and denoted by
$\Sigma^{i_c}$. Again, passing to extrinsic complexification, 
from the relation $(\Sigma^{i_c})^c=\Sigma^{i_c}\times \C_\zeta^n$ and
using the principle of analytic continuation for complex algebraic 
subsets, it follows that $\Sigma^{i_c}$ is irreducible. 

\subhead 9.6.~Locus of CR points \endsubhead
We denote by $\Sigma_{CR}$ the set of points $z\in\Sigma_{reg}$ at
which the rank of the $\sigma\times n$ complex matrix $({\partial P_i\over
\partial z_k}(z,\bar z))_{1\leq i\leq \sigma, \, 1\leq k\leq n}$ 
is maximal, 
hence locally constant. It is called the set of {\it CR points} of
$\Sigma$. Of course, its extrinsic complexification $(\Sigma_{CR})^c$,
the set of points $(z,\zeta)\in \Sigma_{reg}^c$ at which the complex
matrix $({\partial P_i\over \partial z_k}(z,\zeta))_{1\leq k\leq n, 
\, 1\leq i\leq \sigma}$ is maximal and locally constant, is also connected.
Let $d_1$ be this generic point. We have $0\leq d_1\leq d$. Notice that 
$\Sigma$ is CR-generic at a Zariski-generic rank if and only if
$d_1=d$. Let the CR-codimension $\codim_{CR} \, \Sigma$ be $d-d_1$.
It is well known that in a neighborhood of a point $p\in \Sigma_{CR}$,
then $\Sigma$ is contained in a holomorphic algebraic submanifold
of codimension $(d-d_1)$, which is the {\it local intrinsic complexification}
$\Sigma_{loc}^{i_c}$ of $\Sigma$. Clearly, $\Sigma_{loc}^{ic}$ is contained
in a unique irreducible complex algebraic subset of $\C^n$. Passing to 
extrinsic complexification to use the connectedness of $\Sigma^c$, one
can check that this irreducible complex algebraic subset does not
depend on $p\in\Sigma_{CR}$ and is nothing else than the (global)
intrinsic complexification $\Sigma^{i_c}$ defined in \S9.5 above.

\subhead 9.7.~Nonlocality of algebraic degeneracy \endsubhead
At an arbitrary point $p\in\Sigma_{CR}$, we define the integer 
$\kappa_{\Sigma,p}$ to be the maximal number of $(1,0)$ vector fields 
$L_1,\ldots,L_j$ with algebraic coefficients in a neighborhood of 
$p$ which induce a complex algebraic foliation of $\Cal V_{\C^n}(p)$  
(possibly with singularities\,; of course, algebraic foliation 
implies algebraicity of the coefficients) such that
\roster
\item"{\bf (1)}"
$L_1,\ldots,L_j$ are tangent to $\Sigma_{CR}\cap \Cal V_{\C^n}(p)$.
\item"{\bf (2)}"
$L_1,\ldots,L_j$ are linearly independent\,:
if $b_1,\ldots,b_j$ are algebraic holomorphic in $\Cal V_{\C^n}(p)$
and if $b_1L_1+\ldots + b_jL_j=0$, then 
$b_1=\ldots=b_j=0$.
\endroster

\noindent
Then the vectors $L_1(q),\ldots,L_{\kappa_{\Sigma,p}}(q)$ are
$\C$-linearly independent for $q$ varying in a Zariski open subset 
of $\Sigma_{CR} \cap \Cal V_{\C^n}(p)$. Evidently, the function 
$\Sigma_{CR}\ni p\mapsto \kappa_{\Sigma,p}\in\N$ is lower 
semi-continuous. 

\proclaim{Lemma~9.8}
Let $p,q\in\Sigma_{CR}$. Then $\kappa_{\Sigma,p}=\kappa_{\Sigma,q}$.
\endproclaim

\demo{Proof}
Let $L_1,\ldots,L_{\kappa_{\Sigma,p}}$ satisfy {\bf (1)} and 
{\bf (2)} above in $\Cal V_{\C^n}(p)$. The coefficients of 
$L_j:=\sum_{k=1}^n a_{j,k}(z) \, {\partial \over \partial z_k}$
being holomorphic algebraic in a neighborhood of $p$, $j=1,\ldots,
\kappa_{\Sigma,p}$, there exists a complex algebraic subvariety 
$E$ of $\C^n$ with $p\not\in E$, such that
\roster
\item"{\bf (3)}"
The coefficients $a_{j,k}(z)$ extend holomorphically along any continuous
path $\gamma$ with origin $p$ which is contained in $\C^n\backslash E$ 
and with endpoint an arbitrary point $q\in \C^n\backslash E$. 
\item"{\bf (4)}"
The extended vector fields at $q$ have holomorphic algebraic 
coefficients, they induce an algebraic foliation and they are linearly 
independent at $q$ in the sense of {\bf (2)} above.
\endroster
Notice that the vector fields
$(L_1)^c,\ldots,(L_{\kappa_{\Sigma,p}})^c$ are tangent to
$(\Sigma_{CR})^c$ in a neighborhood of $p$. Let $q\in\Sigma_{CR}$ be
arbitrary. Remember that if $A$, $B$ are complex algebraic sets with
$A$ irreducible, $A_{reg}\backslash B$ is connected and nonempty unless $B$
contains $A$.  Hence $(\Sigma_{CR})^c\backslash (E\times \C_\zeta^n)$
is connected. By {\bf (3)}, there exists a continous path $\gamma$
from $(p,\bar p)$ to $(q,\bar q)$ running in
$(\Sigma_{CR})^c\backslash (E\times \C_\zeta^n)$.  By the principle of
analytic continuation, it follows that the extended vector fields
along $\gamma$, say $L_1',\ldots,L_{\kappa_{\Sigma,p}}'$ are tangent
to $(\Sigma_{CR})^c\cap \Cal V_{\C^{2n}}(q,\bar q)$. We deduce
$\kappa_{\Sigma,p}\leq \kappa_{\Sigma,q}$, and then
$\kappa_{\Sigma,p}= \kappa_{\Sigma,q}$, as desired.
\qed
\enddemo

\subhead 9.9.~Local straightening property \endsubhead
Let $\kappa_{\Sigma}$ denote the common value of all the
$\kappa_{\Sigma,p}$ for $p\in\Sigma_{CR}$ and call it the
{\it degree of algebraic degeneracy} of $\Sigma$. 
The following statement achieves to explain Theorem~2.5.

\proclaim{Theorem~9.10}
There exists a proper real algebraic subvariety $F$ of $\Sigma$ with
$\Sigma\backslash F\subset \Sigma_{CR}$ such that in a neighborhood of
every point $p\in\Sigma\backslash F$, there exist local holomorphic
algebraic coordinates in which $\Sigma$ is of the form
$\Delta^{\kappa_\Sigma}\times\underline{\Sigma}_p$, where
$\underline{\Sigma}_p\subset \C^{n-\kappa_\Sigma}$ is a piece of smooth CR 
real algebraic subset satisfying $\kappa_{\underline{\Sigma}_p}=0$.
\endproclaim

\demo{Proof}
We check first $\kappa_{\underline{\Sigma}_p}=0$.
Otherwise, let $\underline{L}$ be a $(1,0)$ vector field
in $\C^{n-\kappa_\Sigma}$ inducing an algebraic foliation tangent
to $0\times \underline{\Sigma}_p$. Of course, 
$\underline{L}$ is linearly independent with 
$L_1=\partial/ \partial z_1,\ldots,L_{\kappa_\Sigma}=\partial
/\partial z_{\kappa_{\Sigma}}$. This contradicts $\kappa_{\Sigma,p}=
\kappa_\Sigma$.

Let us now establish the product property. According to the preceding
considerations, there exist multivalued global vector fields
$L_1,\ldots,L_{\kappa_\Sigma}$ defined outside a complex algebraic
subset $E$ of $\C^n$. Enlarging $E$ if necessary, we can assume that
the vectors $L_1(q),\ldots,L_{\kappa_\Sigma}(q)$ are linearly
independent at every point $q\in\C^n\backslash E$. Let
$\Sigma_{NCR}\supset
\Sigma_{sing}$ be the set of non-CR points of $\Sigma$. 
Clearly, $\Sigma_{NCR}$ is a proper real algebraic subvariety of $\Sigma$.
We set 
$F:= \Sigma_{NCR}\cup (E\cap \Sigma)$.  According to \S9.6 above,
$\Sigma_{CR}$ is CR-generic in its local intrinsic complexification
$(\Sigma_{CR} \cap \Cal V_{\C^n}(p))^{i_c}=
\Sigma^{i_c}\cap \Cal V_{\C^n}(p)$, which is a smooth complex algebraic
subvariety of codimension equal to $\codim_{CR} \Sigma^{i_c}$.  The
vector fields $L_1,\ldots,L_{\kappa_\Sigma}$ being tangent to
$\Sigma_{CR}$, their restriction to $\Sigma^{i_c} \cap \Cal
V_{\C^n}(p)$ is tangent to it. To prove Theorem~9.10, we can therefore
reason directly inside $\Sigma^{i_c} \cap \Cal V_{\C^n}(p)$. In other
words, we come to the following statement.
\enddemo

\proclaim{Lemma~9.11}
Let $M\subset \C^n$ be a \text{\rm CR-generic} real algebraic submanifold, 
let $p\in M$ and let $L_1,\ldots,L_\kappa$, $\kappa\geq 1$, be
$(1,0)$ vector fields with algebraic coefficients such that
\roster
\item"{\bf (1)}"
The vectors $L_1(q),\ldots,L_\kappa(q)$ are linearly independent.
\item"{\bf (2)}"
The induced foliations $\Cal F_{L_1},\ldots, \Cal F_{L_\kappa}$ are 
algebraic in $\Cal V_{\C^n}(p)$.
\item"{\bf (3)}"
The vector fields $L_1,\ldots, L_\kappa$ are tangent to $M\cap 
\Cal V_{\C^n}(p)$.
\endroster
Then there exist local algebraic coordinates at $p$ in which $M=
\Delta^\kappa\times \underline{M}$, where $\underline{M}\subset 
\C^{n-\kappa}$ is a CR-generic real algebraic submanifold.
\endproclaim

\demo{Proof}
According to \S6.1, there exist local coordinates $t=(w,z)\in\C^m\times
\C^d$ in which $p=0$ and $M\: z_l=\bar Q_l(w,\bar t)$, 
$l=1,\ldots, d$. We shall identify $(t_1,\ldots,t_m)=(w_1,\ldots,w_m)$.
Of course $\kappa\leq m$. Let $L_1=\sum_{k=1}^n a_k(t) \, \partial
/\partial t_k$. We can assume that $a_1(0)\neq 0$. Since $a_k/a_1$ is
still holomorphic algebraic, $k=2\ldots,n$, we come down to 
$L_1=\partial / \partial w_1+\sum_{k=2}^n a_k(t) \,
\partial / \partial t_k$. Due to $a_1(t)\equiv 1$, 
the following lemma puts in concrete form the
assumption of algebraicity of $\Cal F_{L_1}$.

\proclaim{Lemma~9.12}
The foliation $\Cal F_{L_1}$ is algebraic if and only if the 
flow of $\Cal F_{L_1}$ is algebraic.
\endproclaim

\demo{Proof} 
For small $u\in\C$ and $t\in\C^n$,
let $\varphi_1(u,t):= \exp (u\, L_1)(t)$, be the local flow of 
$L_1$ in a neighborhood of $0$. If it is algebraic, then the algebraic
biholomorphism $\Phi(w_1,t_2,\ldots,t_n):= \varphi_1(w_1,0,t_2,\ldots,
t_n)$ produces some new straightened algebraic coordinates in 
which the integral curves of $L_1$ are the affine lines
$\{t_2=c_2,\ldots,t_n=c_n\}$, since $\Phi_*(\partial / \partial w_1)
=L_1$. This shows the algebraicity of $\Cal F_{L_1}$.

Conversely, let $\Phi\: \C\times \C^{n-1}\ni (x,y)\mapsto
\Phi(x,y) \in\C^n$, $\Phi(0)=0$, be a local algebraic coordinate
system in which the lines $\{y=c\}$ coincide with the 
complex integral curves of $L_1$.
As $a_1(t)\equiv 1$, for fixed $u_*$, the flow of $L_1$ at time
$u_*$, restricted to the algebraic hypersurface $\{w_1=0\}$, 
namely the map $(t_2,\ldots, t_n)\mapsto \varphi_1(u_*,0,
t_2,\ldots,t_n)$, sends $\{w_1=0\}$ onto the algebraic hypersurface
$\{w_1=u_*\}$. Thanks to the coordinates defined by $\Phi$,
we deduce that
\roster
\item"{\bf 1.}"
The maps $(t_2,\ldots,t_n)\mapsto \varphi_1(w_{1*}, 0,t_2,\ldots,t_n)$
are algebraic.
\item"{\bf 2.}"
The maps $w_1\mapsto \varphi_1(w_1,0,t_{2*},\ldots,t_{n*})$ 
are algebraic.
\endroster
According to the separate algebraicity principle [3], the map 
$(w_1,t_2,\ldots,t_n)\mapsto \varphi_1(w_1,0,t_2,\ldots,t_n)$ is then
algebraic. Thanks to $a_1(t)\equiv 1$, we have
$\varphi_1(w_1,t_1,t_2,\ldots$ $t_n)\equiv 
\varphi_1(w_1+t_1,0,t_2,\ldots, t_n)$, so the total 
flow of $L_1$ is algebraic, as desired.
\qed
\enddemo

\demo{End of proof of Lemma~9.11}
Since the straightened field $\partial / \partial w_1$ is tangent to
$M$, in coordinates like in \S6.1, the functions $\bar Q_l$ are
independent of $(w_1,\bar w_1)$. Thus $M$ is the product $\Delta\times
\underline{M}$, where $\underline{M}\subset \C^{n-1}$ is given by
$z_l=\bar Q_l(w_2,\ldots,w_m, \bar t_2,\ldots,\bar t_n)$, $l=1,\dots,
d$.  If $\kappa\geq 2$, let $L_2=\sum_{k=1}^n a_k(t) \, \partial
/\partial t_k$ be linearly independent with $\partial / \partial w_1$
at $0$. We can assume that $a_2(t)\equiv 1$, namely
$L_2=a_1(t) \, \partial / \partial w_1+\partial / \partial w_2+
\sum_{k=3}^n a_k (t) \, \partial / \partial t_k$.
Thanks to Lemma~9.12 above, its flow $\varphi_2(u,t):= \exp (uL_2)(t)$
is algebraic. Then $\Phi(t_2,t_3,\ldots,t_n):=
\exp (t_2L_2)(0,0,t_3,\ldots,t_n)$ is algebraic. Denote 
$\Phi:=(\Phi_1,\Phi_2,\ldots$ $\Phi_n)$ and
$\underline{\Phi}:=(\Phi_2,\ldots,\Phi_n)$. We have $\hbox{det} 
({\partial \underline{\Phi}_j\over\partial t_k}(0) 
)_{2\leq j,k\leq n}\neq 0$. The map $\underline{\Phi}$ therefore 
induces a local algebraic
biholomorphism of $0\times \C^{n-1}$ fixing $0$.  Since $L_2$ is tangent
to $M$, we have $(0,0,t_3,\ldots,t_n)\in M$ iff
$\Phi(t_2,t_3,\ldots,t_n)\in M$. 
As the equations of $M$ do not depend on $(w_1,\bar
w_1)$, it follows that $\underline{\Phi}(0,t_3,\ldots, t_n)\in
\underline{M}$ iff $\underline{\Phi} (t_2,t_3,\ldots,t_n)\in
\underline{M}$.  Consequently the vector field $\underline{L}_2:=
\underline{\Phi}_*(\partial / \partial t_2)$ is tangent 
to $\underline{M}$. Further, the foliation $\Cal F_{\underline{L}_2}$ is
algebraic, since $\underline{\Phi}$ is. In the coordinates over
$\C^{n-1}$ defined by $\underline{\Phi}$, we have
$\underline{L}_2=\partial / \partial t_2$. As above, we deduce that
the equations of $\underline{M}$ are independent of $(w_2,\bar w_2)$,
whence $M=\Delta^2\times
\underline{\underline{M}}$, with $\underline{\underline{M}}\subset
\C^{n-2}$ being CR-generic. If $\kappa\geq 3$, we proceed analogously, 
{\it etc.} up to $L_\kappa$, which completes the proofs of Lemma~9.11
and of Theorem~9.10.
\qed
\enddemo
\enddemo

\subhead 9.13. Holomorphic degeneracy \endsubhead
Let $\Sigma$ be an irreducible real algebraic subset of $\C^n$.  Let
$p\in \Sigma_{CR}$. 
Let $\kappa_{\Sigma,p}^{hol}$ denote the maximal
number of $(1,0)$ vector fields $L_1,\ldots,L_j$ with holomorphic
coefficients which are tangent to $\Sigma_{CR}\cap
\Cal V_{\C^n}(p)$ and linearly independent over $\Cal H(\Cal V_{\C^n}(p))$
({\it cf.}~[2,15]). Clearly, we have $\kappa_{\Sigma,p}^{hol}\geq
\kappa_\Sigma$, since we do not require the induced foliations to be
algebraic. Nevertheless, we can deduce from Theorem~2.5 that these
two integers coincide.

\proclaim{Corollary~9.14}
For all $p\in\Sigma_{CR}$, we have $\kappa_{\Sigma,p}^{hol}=
\kappa_{\Sigma,p}=\kappa_\Sigma$.
\endproclaim

\demo{Proof}
Set $\chi:=\kappa_{\Sigma,p}^{hol}$.  Let us choose $q\in\Cal
V_{\C^n}(p)$ with the vectors $L_1(q), \ldots, L_{\chi}(q)$ being linearly
independent, in a neighborhood of which $\Sigma=
\Delta^{\kappa_\Sigma}\times \underline{\Sigma}_q$.
As $\chi> \kappa_\Sigma$, a linear combination of the $L_j$'s is
tangent to $0\times \underline{\Sigma}_q$.  We are thus reduced to
prove the corollary in the case where $\kappa_\Sigma=0$ and $\chi\geq 1$. By
contradiction, assume that $\kappa_\Sigma=0$ but there exists a
holomorphic vector field $L$ with $L(p)\neq 0$ tangent to $\Sigma\cap
\Cal V_{\C^n}(p)$.  Let $\varphi(u,t):=\exp (uL)(t)$ be the flow of
$L$, written in coordinates $t$ vanishing at $p$. Set $u_{a,b}(t):=
b+\sum_{k=1}^n a_k \, t_k$, with small $a\in\C^n$ and $b\in\C$. 
According to the separate
algebraicity principle, $\varphi(u,t)$ is algebraic if and only if
$t\mapsto \varphi(u_{a,b}(t),t)$ is algebraic for all $(a,b)$.  If
$\varphi(u,t)$ was algebraic, there would exist an algebraic
biholomorphism straightening $L$ in $\partial / \partial t_1$ ({\it
see} the proof of Lemma~9.12), whence $\kappa_\Sigma\geq 1$, which is
untrue. Thus there exist $a$ and $b$ such that $t\mapsto \varphi
(u_{a,b}(t),t)$ is a nonalgebraic biholomorphic self-map of
$(\Sigma,p)$. Then Corollary~3.3 implies $\kappa_\Sigma \geq 1$,
contradiction.
\qed
\enddemo

\subhead 9.15.~Link with the defining equations \endsubhead
Finally, we describe a useful mean of calculating $\kappa_\Sigma$
({\it cf.}~[2,15]). Without loss of generality, we can pick a small
connected piece $M$ of $\Sigma_{CR}$ centered at one of its points $p$
at which $M$ is CR-generic and given by holomorphic algebraic
equations $\bar z_l=Q_l(\bar w,t)$, $l=1,\dots, d$. Let us write
$Q_l(\bar w,t)=\sum_{\beta\in\N^m} \bar w^\beta \, Q_{l,\beta}
(t)$. Let $\Cal Q(t)$ denote the $d\times \infty$ matrix of
holomorphic algebraic functions $(Q_{l,\beta}(t))_{1\leq l\leq d, \,
\beta\in\N^m}$. After a dilatation of coordinates, we can assume that
all the $Q_{l,\beta}$ are holomorphic in $\Delta^n$. Let us study
abstractly the situation where we are given such a denumerable
collection $\Phi(t):= (\varphi_k(t))_{k\in\N}\in\C^\infty$ of
holomorphic functions over $\Delta^n$. By the Jacobian $\hbox{Jac} \,
\Phi(t)$ of $\Phi(t)$, we understand the $n\times \infty$ matrix
$({\partial \varphi_k\over \partial t_l}(t))_{1\leq l\leq n, \, k\in
\N}$. Let $i$ be an integer with $1\leq i\leq n$. We consider the
denumerable collection of all $i\times i$ minors of $\hbox{Jac}
\, \Phi(t)$. Let $\chi$ denote the generic rank of the map $\Phi(t)$. 
In other words, all the $(\chi+1)\times (\chi+1)$ minors of
$\hbox{Jac} \, \Phi(t)$ vanish identically over $\Delta^n$ and 
there exists a $\chi\times \chi$ minor of $\hbox{Jac} \, \Phi(t)$ 
that does not vanish identically. We shall denote by $\text{\rm gen-rk}_\C
\, \Phi=\chi$ the generic rank of $\Phi$. The set $E\subset \Delta^n$ of
points where all the $\chi\times \chi$ minors vanish is then a proper
real algebraic subset of $\Delta^n$. According to the constant rank 
theorem (also valuable for a denumerable collection of holomorphic 
functions), for each point $p\in \Delta^n\backslash E$, there exists 
a neighborhood $V$ of $p$ contained in $\Delta^n\backslash E$ such that
for each $q\in V$, the set $\Cal F_q:= \{t\in V\: \varphi_k(t)=
\varphi_k(q) \ \forall k\in \N\}$ is a complex algebraic variety and
the union of the $\Cal F_q$ equips $V$ with a holomorphic algebraic
foliation of dimension $n-\chi$. Applying this to $\Cal Q(t)$,
we deduce the following fundamental result.

\proclaim{Theorem~9.16}
Let $M$ be a connected local piece of a real algebraic CR-generic
submanifold of $\C^n$ of CR dimension $m$ and of codimension $d$
passing through the origin which is given in coordinates
$t=(w,z)\in\C^m\times \C^d$ by the $d$ scalar equations $\bar
z_l=Q_l(\bar w,t)=\sum_{\beta\in\N^m}
\bar w^\beta \, Q_{l,\beta}(t)$, $l=1,\ldots,d$, let 
$\Cal Q(t)$ denote the $d\times \infty$ matrix of holomorphic
algebraic functions $(Q_{l,\beta}(t))_{1\leq l\leq d, \,
\beta\in\N^m}$. Let $\chi_M$ denote the generic rank of 
$\Cal Q(t)$. Then the degree of algebraic degeneracy $\kappa_M$ of $M$
defined in \S9.7 above has the following properties\,:
\roster
\item"{\bf (1)}"
The relations $\kappa_M+\chi_M=n$ and $\kappa_M\leq m$ hold.
\item"{\bf (2)}"
The set of points $t\in M$ where the rank of $\Cal Q(t)$ is $< \chi_M$
is a proper real algebraic subvariety of $M$.
\endroster
\endproclaim

\demo{Proof}
It remains only to check that $\kappa_M=n-\chi_M$. At first, applying
Theorem~9.10, we see immediately that $\chi_M\leq
n-\kappa_M$. Conversely, at a point $p$ close to the origin at which
the matrix $\Cal Q(t)$ is of constant rank, a neighborhood $V$ of $p$
algebraically foliates by $(n-\chi_M)$-dimensional leaves defined by
$\Cal F_q=\{t\in V\: Q_{l,\beta}(t)=Q_{l,\beta}(q), \ \forall \, l,
\beta\}$. Of course, such a foliation is algebraically biholomorphic
to a product by a polydisc $\Delta^{n-\chi_M}$. It remains only to
show that the foliation is {\it tangent} to $M$, namely that one of its
leaf is entirely contained in $M$ if and only if it intersects
$M$. But suppose $q\in M$, {\it i.e.} $z_q=\bar Q(w_q,\bar t_q)$, and let
$t\in \Cal F_q$ which by definition satisfies $\bar Q_{l,\beta}(\bar
t)=\bar Q_{l,\beta}(\bar q)$ for all $l$, $\beta$. We deduce $z_q=\bar
Q(w_q, \bar t)$. Remember that there exists an invertible $d\times d$
matrix $a(t,\tau)$ satisfying $z-\bar Q(w,\tau)\equiv a(t,\tau) \,
(\xi-Q(\zeta, t))$. We deduce $\bar z=Q(\bar w,t_q)$.  Finally, $\bar
z=Q(\bar w, t)$, again since $Q_{l,\beta}(t)=Q_{l,\beta}(q)$ for all
$l$, $\beta$. Hence $t\in M$. This shows that $\kappa_M\geq n-\chi_M$, 
as desired.
\qed
\enddemo

\head \S10. Segre-transversality \endhead 

\subhead 10.1.~Definitions\endsubhead Our goal in this paragraph
is to establish that up to and including codimension two, 
Segre-transversality is equivalent to minimality, at Zariski-generic
points. At first, we need preliminary definition. As above, let $M$ be a
CR-generic submanifold of $\C^n$ with $m:=\dim_{CR} \, M$ and
$d=\codim_{\R} \, M$. We assume that $M$ is real analytic or real 
algebraic. Let $p\in M$. Let $S_{\bar q}$
denote the Segre varieties associated with 
$q\in\Cal V_{\C^n}(p)$, {\it cf.}~\thetag{5.3}. 
Remember that $q\in S_{\bar p}$ iff $p\in S_{\bar q}$.
Following [5] and slightly refining the
notions, we shall say that $M$ is called Segre-transversal {\it at}
$p$ if for each neighborhood $V={\Cal V}_{\C^n}(p)$ of $p$, $\exists
\, k \in \N_*$, $\exists \, p_1,\ldots,p_k\in S_{\bar{p}} \cap V$ such
that
$$ 
T_p S_{\bar{p}_1} +
\cdots +T_p S_{\bar{p}_k} =T_p \C^n. 
\tag 10.2
$$ 
Also, $M$ is called Segre-transversal {\it in} $p$ if for each
neighborhood $V={\Cal V}_{\C^n}(p)$ of $p$, $\exists \, k \in \N_*$,
$\exists \, q\in V$, $\exists \, p_1,\ldots,p_k\in S_{\bar{q}} \cap V$
such that
$$ 
T_q S_{\bar{p}_1} + \cdots +T_q S_{\bar{p}_k} =T_q
\C^n. 
\tag 10.3
$$

\remark{Remarks}
{\bf 1.}~Obviously, Segre-transversality {\it at} $p$ implies 
Segre-transversality {\it in} $p$.

{\bf 2.}~We observe that Segre-transversality at $p$ {\it does not
entail} that the second {\it Segre set} ({\it cf.} [2])
$S_{\bar{p}}^2:=\cup_{q\in
S_{\bar{p} \cap V}} S_{\bar{q}} \cap V$ contains an open set in
$\C^n$, as shows the following example in $\C^3$, $M\:
z_1=\bar{z}_1+iw\bar{w}$, $z_2=\bar{z}_2+iw^2\bar{w}^2$, or, more 
generally, as shows any $M$ of codimension $d>2m$.

{\bf 3.}~Nevertheless, it is known and clear 
that {\it Segre-transversality at
$p$ entails minimality at $p$}, because it is easy to 
check in coordinates that all the Segre varieties
$S_{\bar{p_1}},\ldots, S_{\bar{p_k}}$ must be contained in the
intrinsic complexification ${\Cal O}_p^{i_c}$ of the CR orbit of $p$
(see also [2,9,14]), whence condition~\thetag{10.2} forces 
${\Cal O}_p^{i_c}$ to be of dimension $n$.
\endremark

\proclaim{Proposition~10.4}
The following properties hold\,:
\roster
\item"{\bf (1)}"
A connected CR-generic manifold $M\subset \C^n$ of codimension
$1$ or $2$ is minimal at a Zariski-generic point if and only if it is 
Segre-transversal in at least one point $p\in M$.
\item"{\bf (2)}"
The CR-generic manifold $M\subset \C^4$ of codimension $3$ given by
$z_1=\bar{z_1}+iw\bar{w}$, $z_2=\bar{z}_2+ iw\bar{w}(w^2+\bar{w}^2)$
and $z_3=\bar{z}_3+iw^3\bar{w}^3$ is minimal at $0$ but not 
Segre-transversal
in any point.
\endroster
\endproclaim

\demo{Proof}
The example in part {\bf (2)} is borrowed from [5]. Part {\bf (1)} is
new for the codimension $d=2$. 
In order to prove this proposition, we plan now to give
characterizations of Segre-transversalities. In fact, a characterization of
Segre-transversality at a point is given in [5], but for $M$ {\it
rigid algebraic} only. We seek a general formulation. Thus, we assume
that $p=0$, $V=\Delta^n$ and that $M$ is given by~\thetag{6.3}.
Let $q\in
\Delta^n$, let $k\in \N^*$, let $p_1,\ldots,p_k\in S_{\bar{q}}\cap
\Delta^n$, {\it i.e.} 
$$
z_{p_j}=\bar{z}_q+i\bar{\Theta}(w_{p_j},\bar{w}_q,\bar{z}_q), \ \ \ \ \
j=1,\ldots, k.
\tag 10.5
$$ 
To compute $T_qS_{\bar{p}_j}$ is easy\,: $T_q
S_{\bar{p}_j}$ is in fact generated by the $m$-vector field
({\it cf.}~\thetag{6.2})
$$ 
{\Cal L}=\frac{\partial }{\partial
w}+i\bar{\Theta}_w( w_q,\bar{w}_{p_j},
z_q-i\Theta(\bar w_{p_j},w_q,z_q)
) \frac{\partial
}{\partial z},
\tag 10.6
$$ 
{\it i.e.} by a collection
$v^l=v^l(w_q,z_q,\bar{w}_{p_j})$ of vectors $l=1,\ldots, m$ of the
form 
$$
(0,\ldots,1,\ldots,0,v_1^l(w_q,z_q,\bar{w}_{p_j}),\ldots,
v_d^l(w_q,z_q,\bar{w}_{p_j}))\in \C^m\times \C^d,
\tag 10.7
$$ 
with $1$ at the $m$-th place, $j=1,\ldots, k$. Let us denote by
$\Lambda(w_q,z_q,\bar{w}_{p_j})$ the linear space generated by the
$v^l$, for $l=1,\ldots, m$. By definition, $M$ is Segre-transversal in $p$
if and only if there exists $q$ arbitrarily close to $p$ and 
$p_1,\ldots, p_k\in S_{\bar q}$ such that
$$ 
\text{\rm Span} \
(\Lambda(w_q,z_q,\bar{w}_{p_1}),\ldots,
\Lambda(w_q,z_q,\bar{w}_{p_k}))=\C^n. 
\tag 10.8
$$ 
Of course it suffices to choose $k=d+1$ to test (10.8). We can also
complexify~\thetag{10.8} and get
$$ 
\text{\rm Span} \
(\Lambda(w,z,\zeta_1),\ldots,\Lambda(w,z,\zeta_{d+1}))=\C^n.
\tag 10.9
$$ 
Introducing all $n\times n$ minors of the system of $m(d+1)$ vectors
$v^l(w,z,\zeta_j)$, $l=1,\ldots,
m$, $j=1,\ldots, d+1$, denoting by $N$ the number of such
minors, we may build a holomorphic map $\X\: {\Cal
M}_{\natural}^{d+1} \to \C^N$ whose components are
these minors and which satisfies $\X(t,0,\ldots,t,0)\equiv 0$.
We may equip with a projection $\pi_{\natural}\: {\Cal M}_{\natural}^{d+1}
\to \C_t^n$ the set
$$ 
{\Cal
M}_{\natural}^{d+1}=\{ (t_1,\tau_1,\ldots,t_{d+1}, \tau_{d+1}) \in
{\Cal M}^{d+1} \: t_i=t_j, 1\leq i,j\leq d+1\}
\tag 10.10
$$ 
such that (10.9) holds at $(w,z)$ if and only if there exists
$p_{\natural} \in \pi_{\natural}^{-1}(w,z)$ with
$\X(p_{\natural})\neq 0$. Consequently (10.9) holds at $t=(w,z)$ if and
only if 
$$
\X(t,\zeta_1,\ldots,t,\zeta_{d+1})\not\equiv_{\zeta_1,\ldots,\zeta_{d+1}}
0\in \C^N. 
\tag 10.11
$$ 
Now writing 
$$
\X(t,\zeta_1,\ldots,t,\zeta_{d+1})= \sum_{\gamma_{\natural}\in
\N_*^{m(d+1)}} \zeta_{\natural}^{\gamma_{\natural}}
\X_{\gamma_{\natural}}(w,z), \ \ \ \ \ \X_{\gamma_{\natural}}(w,z)\in
\C^N, 
\tag 10.12
$$ we have got that there exist $\gamma_{\natural}\in
\N_*^{m(d+1)}$ and $(w,z)\in \Delta^n$ such that
$\X_{\gamma_{\natural}}(w,z)\neq 0$. As a usual consequence of
analyticity (here of $\X$) we see that (10.9) holds at $(w_q,z_q)$
if and only if (10.9) holds for every $(w,z)\in \Delta^n$ minus a
proper complex analytic subset.
We thus have proved\,:

\proclaim{Proposition 10.13} 
A connected real analytic or real algebraic
CR-generic manifold $M$ is Segre-transversal in a point if and
only if it is Segre-transversal in every point. In that case $M$ is
Segre-transversal at every point outside a proper real analytic 
or real algebraic subvariety. In this case, 
$M$ is be called \text{\rm Segre-transversal}.
\endproclaim 

\demo{Proof of Proposition~10.4 {\bf (1)}} 
If $d=1$, it is very easy to check that $M$ is Segre-transversal if
and only if $M$ is not Levi-flat if and only if it is minimal at a 
Zariski-generic point, {\it
cf.} [5]. Thus, let $d=2$ and assume that $M$ is minimal at a 
Zariski-generic point. To
simplify the notation, let us assume that $m=1$. Then at a
Zariski-generic point $p\in M$, the vector-valued Levi-form of $M$ has
rank $1$. Further, the second order Lie brackets of $T^cM$ must
complete the tangent space to $M$ at a Zariski-generic point
(otherwise, the distribution spanned by $T^cM$ and its first order Lie
brackets is involutive and then $M$ locally foliates in three-dimensional CR
orbits, contradiction). Thus, there are holomorphic
coordinates $(w,z_1,z_2)$ vanishing at such a Zariski-generic point
$p$ in which the equations of $M$ are $z_2=\bar z_2+iw\bar w \,
[1+O(\v w\v) +O(\bar z_2,\bar z_3)]$ and $z_3=\bar z_3+iw\bar w
\, [w+\bar w+O(\v w\v^2)+O(
\bar z_2,\bar z_3)]$. Choosing then a point $q=(w_q,0,0)$ with 
$w_q\neq 0$ small, it is easy to check that condition~\thetag{10.8} is
satisfied for suitable points $p_1,p_2,p_3\in S_{\bar q}$ with
$w_{p_1}=0$, and $w_{p_2}$, $w_{p_3}$ small. The general case $m\geq
2$ is similar or can be reduced to the case $m=1$
by slicing. This completes the proof of Proposition~10.4.
\qed
\enddemo
\enddemo

\head \S11. Algebraicity of the reflection mapping and transcendence
degree of the mapping \endhead

\subhead 11.1.~Introduction\endsubhead
Now, to end-up this article, we come to the description of an
equivalent formulation of our main Theorem~2.5. An interesting
invariant of which to show analyticity (in case $f$ is in ${\Cal
C}_{CR}^{\infty}$ and $M$, $M'$ are real analytic) or algebraicity (in
case $f$ is already holomorphic, and $M$, $M'$ are real algebraic), is
the so-called ``{\it reflection mapping}'' ({\it see} [11], where the
interest of its study {\it with no nondegeneracy condition on $M'$}
was pointed out for the first time). The reflection mapping is the
holomorphic map $\C^n\times \C^{n'}\ni (t,\tau') \mapsto
\xi'-Q'(\zeta',f(t)) \in\C^{d'}$, where $\bar z'=Q'(\bar w',t')$ 
is a local system of $d'$ defining equations for $M'$ as in \S6.1 
above. Let us denote it by $\Cal R_f'(t,\tau')$, where $\tau'=
(\zeta',\xi')$. Clearly, $\Cal R_f'$ is relatively algebraic 
with respect to $\tau'$ and {\it a priori} only holomorphic with 
respect to $t$.

\proclaim{Theorem~11.2} \text{\rm ($d=1$\,:
[12])} Let $h\: M\to M'$ be a holomorphic map of generic rank $n$
between two small pieces of connected somewhere minimal CR-generic
real algebraic submanifolds of $\C^n$ of the same CR dimension $m\geq
1$ and of codimension $d\geq 1$. Then the reflection mapping $\Cal R_h'$
associated with $h$ and an arbitrary local system of coordinates for
$M'$ as above is complex algebraic with respect to both variables
$(t,\tau')$.
\endproclaim

\demo{Proof}
Using the biholomorphic invariance of Segre varieties, one can check
that the algebraicity of $\Cal R_h'$ is preserved under changes of
algebraic coordinates fixing a point and under small shifts of a
center point $p'\in M'$ (exercise). Thus, after a slight
delocalization, namely after an arbitrarily small shift of the center
point, and after a change of algebraic coordinates as in Theorem~9.10,
we can assume that $f\: (M,p)\to (M',p')$ is a local biholomorphic
mapping with $M'= \Delta^{\kappa'}\times \underline{M}'$,
$\kappa':=\kappa_{M'}$.  Composing $f$ with a projection, we get a 
submersive local holomorphic map
$\underline{f} \: M\to
\underline{M}'$, with $\underline{M}'$ minimal for inclusion. Since  
$\kappa_{\underline{M}'}=0$, Theorem~3.1 yields that $\underline{f}$
is algebraic. Notice that there exist coordinates
$t'=(w',z')=(u',v',z')\in\C^{m'-\kappa'}\times
\C^{\kappa'}\times\C^{d'}$ as in \S6.1
such that the equations of $(M',p')$ are independent of $(v',\bar
v')$, namely $\bar z'=Q'(\bar u',u',z')$. Let
$(u',v',z')^c=:(\mu',\nu',\xi')$.  Finally, the reflection mapping
$\Cal R_f'(t,\tau') =\xi'-Q'(\mu',\underline{f}(t))\in\C^{d'}$ is then
obviously algebraic, because $\underline{f}$ is. The proof of
Theorem~11.2 is complete.
\qed
\enddemo

\subhead 11.3.~Algebraicity of the reflection mapping with no rank
assumption \endsubhead 

\noindent
More generally, by similar arguments, we shall now establish
that {\it the algebraicity of the reflection mapping associated with
$f\: M \to M'$ is in fact equivalent to Theorem~2.5.} An important
assumption will be that $M'$ is here minimal for inclusion containing
$f(M)$.

\proclaim{Theorem~11.4} 
Let $f$ be a local holomorphic map $(M,p)\to (M',p')$ between two real
algebraic CR-generic manifolds with $(M,p)$ minimal at a
Zariski-generic point, with $(M',p')$
\text{\rm minimal for inclusion} containing $f(M,p)$ and given by
$d'$ equations of the form $\bar z_{l}'=Q_{l}'(\bar w',t')$,
$l=1,\ldots,d'$. Then the reflection mapping\,:
$$
(\C^n,p)\times (\C^{n'},\bar p') \ni (t,\bar{\nu}') \mapsto
(\xi_{l}'-Q_{l}'(\zeta',h(t)))_{1\leq l\leq d'} \in
\C^{d'}
\tag 11.5
$$
extends as an algebraic map of $(t,\tau') \in \C^n \times
\C^{n'}$. Conversely, if the above-defined reflection mapping is
algebraic, then $\nabla^{tr}(f)\leq \kappa_{M'}$.
\endproclaim

\demo{Proof} 
Since $M'$ is minimal for inclusion and since the exceptional set of
points $q'$ in a neighborhood of which $M'\cap {\Cal V}_{\C^{n'}}(q')$
is not equivalent to a product $\Delta^{\kappa_{M'}}\times
\underline{M}'$ is a proper real algebraic 
subvariety of $M'$ ({\it cf.}~Theorems~9.10 and~9.16), 
then $f(M)$ encounters in fact the set of points
$q'$ for which $M'\cap {\Cal V}_{\C^{n'}}(q')\cong_{\Cal A}
\Delta^{\kappa_{M'}}\times
\underline{M}'$. As in the proof of 
Theorem~11.2, we then deduce the algebraicity of $\Cal R_f'$.

Conversely, let $M'$ be minimal for inclusion, let $M'=
\Delta^{\kappa'}\times\underline{M}'$ in a neighborhood of $q'$ close
to $p'$, with $\kappa':=\kappa_{M'}$, choose coordinates $(u',v',z')$
vanishing at $q'$ as above and assume that
$\xi'-Q'(\mu',\underline{f}(t))$ is algebraic. Differentiating with
respect to $\mu'$, we deduce that the functions ${1\over \beta!}
\, [\partial_{\mu'}^\beta Q'(\mu',\underline{f}(t))]_{\mu'=0}=:
Q_\beta'(\underline{f}(t))$, $\beta\in\N^{m'-\kappa'}$, are
algebraic. By Theorem~9.16~{\bf (2)}, it follows from the assumption
$\kappa_{\underline{M}'}=0$ that the set $E'$ of points $(u',z')\in
\underline{M}'$ at which all $d'\times d'$ minors of the Jacobian 
matrix of $\Cal Q'(t')=(Q_{l,\beta}'(u',z'))_{1\leq l\leq d', \,
\beta\in\N^{m'-\kappa'}}$ vanish is a
proper real algebraic subvariety of $M'$. By minimality for inclusion
of $M'$, it follows that $\underline{f}(M,p)$ is not contained in
$E'$.  Using then the algebraic implicit function theorem we deduce
that $\underline{f}$ is algebraic. This implies that $\nabla^{tr}(f)
\leq \kappa'$, as desired.
\qed 
\enddemo

\subhead 11.6.~Algebraic approximation \endsubhead
Let us conclude this paper with a nice application of Theorem~11.4.
We need some preliminary. Under the assumptions of Theorem~11.4, let
$t\in\C^n$, $t'\in\C^{n'}$ be coordinates vanishing at $p$, $p'$.  Let
us denote $r(t,\tau):=z-\bar Q(w,\tau)$ and $r'(t',\tau'):= z'-\bar
Q'(w',\tau')$, whence $\bar r(\tau,t)=\xi-Q(\zeta,t)$ and $\bar
r'(\tau',t')=\xi'-Q'(\zeta',t')$. Because $M$ and $M'$ are real, there
exist a $d\times d$ matrix $a(t,\tau)$ of holomorphic algebraic
functions with $a(0)=-I_{d\times d}$ and $r(t,\tau)\equiv a(t,\tau) \,
\bar r(\tau,t)$, and similarly $a'(t',\tau')$ with 
$a'(0)=-I_{d'\times d'}$ and 
$r'(t',\tau')\equiv a'(t',\tau') \, \bar r'(\tau',t')$. Since $(M,p)$
is mapped in $(M',p')$, there exist a $d'\times d$ matrix $b(t,\tau)$
of holomorphic functions such that $r'(f(t),\bar f(\tau))\equiv
b(t,\tau) \, r(t,\tau)$. Here, $f(t)\in\C\{t\}^{n'}$ is a power series
vanishing at $0$. Finally, we set
$Q'(\zeta',t'):=\sum_{\beta\in\N^{m'}} {\zeta'}^\beta \,
Q_\beta'(t')$. Then it follows from Theorem~11.4 that for each
$\beta\in\N^{m'}$, the function $\varphi_\beta'(t):=
Q_\beta'(f(t))={1\over \beta!} [\partial_{\zeta'}^\beta Q'(\zeta',
f(t))]_{\zeta'=0}$ is algebraic. By applying the algebraic
approximation theorem of Artin [1] to the collection of holomorphic
algebraic equations $Q_\beta'(t')-\varphi_\beta'(t)=0$, we deduce that
for each integer $N\in\N_*$, there exists a holomorphic algebraic
mapping $F\: (\C^n,p)\to (\C^{n'},p')$, namely $F(t)\in\C\{t\}^{n'}$,
which is algebraic, such that $Q_\beta'(F(t))=\varphi_\beta'(t)$ for
all $\beta\in\N^{m'}$ and $F(t)\equiv f(t) \ (\hbox{mod} \, \v
t\v^N)$.  We observe that it follows that $F$ maps $(M,p)$ in
$(M',p')$.

\proclaim{Lemma~11.7}
Suppose that $F(t)\in\C\{t\}^{n'}$ satisfies $Q_\beta'(F(t))=
Q_\beta'(f(t))$. Then $r'(F(t),\bar F(\zeta, Q(\zeta,t)))\equiv 0$ in
$\C\{\zeta,t\}^{n'}$.
\endproclaim

\demo{Proof}
Indeed, if we denote $f=(g,h)\in\C^{m'}\times\C^{d'}$, we have first
$\bar h(\zeta,Q(\zeta,t))\equiv \sum_{\beta\in\N^{m'}}
\bar g(\zeta,Q(\zeta,t))^\beta \, Q_\beta'(F(t))$. Let us write 
$F=(G,H)\in\C^{m'}\times\C^{d'}$. First, from the relation
$r'(F(t),\bar f(\tau))\equiv a'(F(t),\bar f(\tau)) \, \bar r'(
\bar f(\tau),F(t))$, we deduce $H(t)\equiv \sum_{\beta\in\N^{m'}}
G(t)^\beta \, \bar Q_\beta'(\bar f(\tau))$. As the coordinates
$(\zeta,t)$ and $(w,\tau)$ are equivalent on $\Cal M$, we deduce
$H(w,\bar Q(w,\tau))\equiv\sum_{\beta\in\N^{m'}}
G(w,\bar Q(w,\tau))^\beta \, \bar Q_\beta'(\bar f(\tau))$.
Finally, we get $H(w,\bar Q(w,\tau))\equiv\sum_{\beta\in\N^{m'}}
G(w,\bar Q(w,\tau))^\beta \, \bar Q_\beta'(\bar F(\tau))$, as
desired.
\qed
\enddemo

\noindent
In particular, we obtain the following interesting application. 
Let $f\: (M,p)\to (M',p')$ be a biholomorphic equivalence between 
two real algebraic CR-generic manifolds such that $(M,p)$ is minimal
at a Zariski-generic point. Then there exists an algebraic holomorphic
equivalence $F\: (M,p)\to (M',p')$. Is such a property also true for
nowhere minimal CR-generic manifolds~?

\enddocument
\end